\documentclass[smallextended,envcountsect]{svjour3}

\usepackage{epsfig} 
\usepackage{amsmath} 
\usepackage{amssymb}  
\usepackage{enumerate}
\usepackage{url}

\newenvironment{pf}{\smallbreak\noindent{\it Proof. }}{\hfill$\Box$\smallbreak}

\usepackage{tikz}
\usepackage{pgfplots}

\newcommand{\id}{\mathrm{Id}}
\newcommand{\overbar}[1]{\mkern 1.5mu\overline{\mkern-1.5mu#1\mkern-1.5mu}\mkern 1.5mu}
\newcommand{\reals}{\mathbb{R}}

\newcommand{\prox}{{\rm{prox}}}
\newcommand{\argmin}{\mathop{\rm argmin}}
\newcommand{\Argmin}{\mathop{\rm Argmin}}

\title{\LARGE \bf
Envelope Functions: Unifications and Further Properties
}

\date{}

\author{Pontus Giselsson \and Mattias F{\"a}lt}

\institute{P. Giselsson (corresponding author) \and M. F{\"a}lt \at
              Department of Automatic Control, Box 118, SE-221 00 Lund, Sweden \\
              Tel.: +46-46-222-\{9744,0847\}, 
              Fax: +46-46-138118\\
              \email{\{pontusg,mattiasf\}@control.lth.se}           
}

\date{Received: date / Accepted: date}

\begin{document}

\maketitle

\begin{abstract}

Recently, the forward-backward and Douglas-Rachford envelope functions were proposed in the literature. The stationary points of these envelope functions have a close relationship with the solutions of the possibly nonsmooth optimization problem to be solved. The envelopes were shown to be smooth and convex under some additional assumptions. Therefore, these envelope functions create powerful bridges between nonsmooth and smooth optimization. 

In this paper, we present a general envelope function that unifies and generalizes these envelope functions. We provide properties of the general envelope function that sharpen corresponding known results for the special cases. We also present an envelope function for the generalized alternating projections method (GAP), named the GAP envelope. It enables for convex feasibility problems with two sets, of which one is affine, to be solved by finding any stationary point of the smooth and under some assumptions convex GAP envelope. 


\end{abstract}


\section{Introduction}

Many convex optimization problems can be solved by finding a fixed-point to a nonexpansive operator. 
This is the basis for many first-order methods such as forward-backward splitting \cite{Combettes2004}, Douglas-Rachford splitting \cite{DouglasRachford,LionsMercier1979}, the alternating direction method of multipliers (ADMM) \cite{Gabay1976,Glowinski1975,BoydDistributed} and its linearized versions \cite{ChambollePock}, the three operator splitting method \cite{Davis_three_splitting}, and generalized alternating projections \cite{GAP_Gubin,GAP_Agmon,GAP_Motzkin,GAP_Eremin,GAP_Bregman} that generalizes \cite{vonNeumann}. 

All these methods seek a fixed-point by performing an averaged iteration of the nonexpansive mapping. 
The averaging is the key to guaranteing convergence of the iterates to a fixed-point of the nonexpansive mapping, see \cite{Combettes2004}. The rate of convergence can, however, be very slow in practice. One way to improve convergence of such methods is to precondition the problem data. This approach has been extensively studied in the literature and has proven very successful in practice; see, e.g., \cite{Benzi_precond,Bramble_Uzawa,Hu_nonlin_Uzawa,GhadimiADMM,gisBoydAut2014metric_select,gisBoydTAC2014metric_select,gisSIAM2015} for a limited selection of such approaches. The underlying idea is to incorporate static second-order information in the respective algorithms.

The performance of the forward-backward and the Douglas-Rachford methods can be further improved by exploiting the properties of the recently proposed forward-backward envelope in \cite{Panos_fb_newton,Panos_fb_quasi_newton} and Douglas-Rachford envelope in \cite{Panos_acc_DR_2014}. As shown in \cite{Panos_fb_newton,Panos_fb_quasi_newton,Panos_acc_DR_2014}, the stationary points of these envelope functions agree with the fixed-points of the corresponding operator. The envelopes are also shown to be convex and to have Lipschitz continuous gradients (under certain assumptions). Therefore, the original nonsmooth problem to be solved using forward-backward splitting or Douglas-Rachford splitting can be solved by finding a stationary point of the corresponding smooth envelope functions. In \cite{Panos_fb_newton,Panos_fb_quasi_newton}, it is shown how truncated Newton methods or quasi-Newton methods can be applied to the forward-backward envelope function to improve local convergence.


A unifying property of forward-backward splitting and Douglas-Rachford splitting (for convex optimization) is that they are averaged iterations of a nonexpansive mapping $S$, where $S=S_2S_1$ is composed of two nonexpansive mappings. These mappings are gradients of functions $f_1$ and $f_2$ respectively, i.e.,  $S_1=\nabla f_1$ and $S_2=\nabla f_2$. What unifies their envelopes is the assumption corresponding to that $f_1$ is twice continuously differentiable. For averaged iteration of such operators, we propose a differentiable envelope function that has the forward-backward and Douglas-Rachford envelopes as special cases. Other special cases include the Moreau envelope and the ADMM envelope (which is a special case of the Douglas-Rachford envelope since ADMM is Douglas-Rachford splitting applied to the Fenchel dual problem, see \cite{Gabay83}).

We analyze this general envelope function in the more restrictive setting of $f_1$ being quadratic, or equivalently $S_1=\nabla f_1$ being affine, i.e., of the form $S_1=P(\cdot)+q$, with $P$ linear. We show that if $P$ is nonsingular, the stationary points of the envelope coincide with the fixed-points of $S=S_2S_1$. We provide quadratic upper and lower bounds to the envelope function that improve corresponding results for the known special cases in the literature. The bounds imply, e.g., that the gradient of the envelope function is always 2-Lipschitz continuous. If in addition the linear operator $P$ that defines $S_1$ is positive semidefinite, the envelope function is convex. Since the fixed-points of $S$ and the stationary points of the envelope coincide, a fixed-point to $S$ can, when $P$ is positive semidefinite, be found by minimizing a smooth and convex envelope function.

In \cite{Panos_fb_newton,Panos_fb_quasi_newton,Panos_acc_DR_2014} it was shown that forward-backward splitting and Douglas-Rachford splitting can be seen as variable metric gradient methods applied to the respective envelope functions. If $S_1$ is affine, they show that it instead is a scaled gradient method with fixed metric. This generalizes also to our setting, i.e., an averaged iteration of a nonexpansive mapping can be interpreted as a scaled gradient method applied to the envelope function. Since the envelope function has nice smoothness properties and is in some cases convex, more efficient methods to find a fixed-point to $S$, or equivalently a stationary point of the envelope, probably exist. For instance, quasi-Newton, nonlinear conjugate gradient, or truncated Newton methods, some of which has been proposed to be used with the forward-backward envelope in \cite{Panos_fb_newton,Panos_fb_quasi_newton} can be used to improve local convergence (see \cite{Nocedal} for details on the methods). Devising new algorithm or suggesting which existing ones that are most efficient is, however, outside the scope of this paper.


We also provide a new envelope function that is a special case of the general envelope, namely the generalized alternating projections (GAP) envelope. Generalized alternating projections \cite{GAP_Gubin,GAP_Agmon,GAP_Motzkin,GAP_Eremin,GAP_Bregman} (which is also referred to as the method of alternating relaxed projections, e.g., in \cite{MARP_Bauschke}) solves feasibility problems involving a finite number of nonempty closed and convex sets. This is done by alternating relaxed projections onto the sets. It can use either under-relaxation, in which the step does not go all the way to the projection point, or over-relaxation when the step goes past the projection point, up towards the reflection point. Our envelope function applies to problems with two sets, with one nonempty closed and convex and one affine. Since the general envelope function always has a Lipschitz continuous gradient, so has the GAP envelope. If in addition, the first relaxed projection (onto the affine set) is an under-relaxation, the GAP envelope is convex. Therefore, all feasibility problems with an affine subspace and a convex set can be solved by minimizing a smooth convex function.


Our contributions are as follows; i) we propose a general envelope function that has several known envelope functions as special cases, ii) we provide properties of the general envelope that sharpen (sometimes considerably) and generalize corresponding known results for the special cases, iii) we provide new insights on the relation between the Douglas-Rachford envelope and the ADMM envelope, iv) we present a new envelope function, the GAP envelope, and characterize its properties.

\section{Preliminaries}

\subsection{Notation}

We denote by $\reals$ the set of real numbers, $\reals^n$ the set of real column-vectors of length $n$, and $\reals^{m\times n}$ the set of real matrices with $m$ rows and $n$ columns. Further $\overbar{\reals}:=\reals\cup\{\infty\}$ denotes the extended real line. We denote inner-products on $\reals^n$ by $\langle\cdot,\cdot\rangle$ and their induced norms by $\|\cdot\|$. We will also use scaled norms $\|x\|_P:=\langle Px,x\rangle$ where $P$ is a positive definite operator (defined in Definition~\ref{def:pos_def}). We will use the same notation for scaled semi-norms, i.e., $\|x\|_P:=\langle Px,x\rangle$ where $P$ is a positive semidefinite operator (defined in Definition~\ref{def:pos_semidef}). The identity operator is denoted by $\id$. The conjugate function is denoted and defined by $f^{*}(y)\triangleq \sup_{x}\left\{\langle y,x\rangle-f(x)\right\}$. The adjoint operator to a linear operator $L~:~\reals^n\to\reals^m$ is defined as the unique operator $L^*~:~\reals^m\to\reals^n$ that satisfies $\langle Lx,y\rangle=\langle x,L^*y\rangle$. The linear operator $L~:~\reals^n\to\reals^n$ is self-adjoint if $L=L^*$. The notation $\argmin_x f(x)$ refers to any element that minimizes $f$ while the notation $\Argmin_x f(x)$ refers to the set of minimizers. Finally, $\iota_C$ denotes the indicator function for the set $C$ that satisfies $\iota_C(x)=0$ if $x\in C$ and $\iota_C(x)=\infty$ if $x\not\in C$.

\subsection{Background}

In this section, we introduce some standard definitions 
that can be found, e.g. in
\cite{bauschkeCVXanal,Rockafellar1997}.

\subsubsection{Operator Properties}
\label{sec:op_prop}
\begin{definition}[Positive semidefiniteness]
A linear operator $L~:~\reals^n\to\reals^n$ is positive semidefinite if it is self-adjoint and all eigenvalues $\lambda_i(L)\geq 0$.
\label{def:pos_semidef}
\end{definition}
\begin{remark}
An equivalent characterization of a positive semidefinite operator is that $\langle Lx,x\rangle\geq 0$ for all $x\in\reals^n$. 
\label{rem:pos_semidef}
\end{remark}
\begin{definition}[Positive definiteness]
A linear operator $L~:~\reals^n\to\reals^n$ is positive definite it is self-adjoint and if all eigenvalues $\lambda_i(L)\geq m$ with $m>0$.
\label{def:pos_def}
\end{definition}
\begin{remark}
An equivalent characterization of a positive definite operator $L$ is that $\langle Lx,x\rangle\geq m\|x\|^2$ for some $m>0$ and all $x\in\reals^n$. 
\label{rem:pos_semidef}
\end{remark}
\begin{definition}[Lipschitz mappings]
  A mapping $T~:~\reals^n\to\reals^n$ is
  $\delta$-\emph{Lipschitz continuous} with $\delta\geq 0$ if
  \begin{align*}
    \|Tx-Ty\|\leq \delta\|x-y\|
  \end{align*}
  holds for all $x,y\in\reals^n$. 
  If $\delta=1$ then $T$ is \emph{nonexpansive} and if $\delta\in[0,1)$
  then $T$ is $\delta$-\emph{contractive}. 
  \label{def:Lipschitz}
\end{definition}
\begin{definition}[Averaged mappings]
  A mapping $T~:~\reals^n\to\reals^n$ is
  $\alpha$-\emph{averaged} if there exists a nonexpansive mapping
  $S~:~\reals^n\to\reals^n$
  and $\alpha\in(0,1]$ such that $T=(1-\alpha)\id+\alpha S$.
\label{def:avg}
\end{definition}
\begin{definition}[Negatively averaged mappings]
  A mapping $T~:~\reals^n\to\reals^n$ is
  $\beta$-\emph{negatively averaged} with $\beta\in(0,1]$ if $-T$ is $\beta$-averaged.
\label{def:negavg}
\end{definition}
\begin{remark}
For notational convenience, we have included $\alpha=1$ and $\beta=1$ in the definitions of (negative) averagedness, which both are equivalent to nonexpansiveness. For values of $\alpha\in(0,1)$ and $\beta\in(0,1)$ averagedness is a stronger property than nonexpansiveness. For more on negatively averaged operators, see \cite{gisSIAM2015} where they were introduced.

Note that if a gradient operator $\nabla f$ is $\alpha$-averaged and $\beta$-negatively averaged. Then it must hold that $\alpha+\beta\geq 1$. This follows immediately from Lemma~\ref{lem:id_sub_grad_mono_avg} and Lemma~\ref{lem:id_sub_grad_mono_negavg} in Appendix~\ref{app:lemmas}.
\end{remark}
\begin{definition}[Cocoercivity]
  A mapping $T~:~\reals^n\to\reals^n$ is $\delta$-cocoercive with
  $\delta > 0$ if $\delta T$ is $\tfrac{1}{2}$-averaged.
  \label{def:cocoercive}
\end{definition}
\begin{remark}
This cocoercivity definition implies that cocoercive mappings $T$ can be expressed as
\begin{align}
T=\tfrac{1}{2\delta}(\id+S)
\label{eq:coco_expression}
\end{align}
for some nonexpansive operator $S$. We also note that 1-cocoercivity is equivalent to $\tfrac{1}{2}$-averagedness (which is also called firm nonexpansiveness).
\label{rem:coco_avg}
\end{remark}
We conclude this subsection with a result relating Lipschitz continuity and cocoercivity to averagedness and negative averagedness.
\begin{proposition}
Suppose that $\nabla f~:~\reals^n\to\reals^n$ is the gradient of some function $f~:~\reals^n\to\reals$. Then the following hold:
\begin{enumerate}[(i)]
\item $\nabla f$ is $\delta$-Lipschitz continuous with $\delta\in[0,1]$ if and only if it is $\tfrac{\delta+1}{2}$-averaged and $\tfrac{\delta+1}{2}$-negatively averaged.
\item $\nabla f$ is $\tfrac{1}{\delta}$-cocoercive with $\delta\in[0,1]$ if and only if it is $\tfrac{1}{2}$-averaged and $\tfrac{\delta+1}{2}$-negatively averaged.
\end{enumerate}
\label{prp:prop_relation}
\end{proposition}
\begin{pf}
Claim {\it{(i)}}: Follows immediately from Lemma~\ref{lem:quad_bound_Lipschitz}, Lemma~\ref{lem:id_sub_grad_mono_avg}, and Lemma~\ref{lem:id_sub_grad_mono_negavg}. Claim {\it{(ii)}}: Lemma~\ref{lem:id_sub_grad_mono_avg}, and Lemma~\ref{lem:id_sub_grad_mono_negavg} imply that $\tfrac{1}{2}$-averagedness and $\tfrac{\delta+1}{2}$-negative averagedness is equivalent to that
\begin{align*}
0\leq f(x)-f(y)-\langle\nabla f(y),x-y\rangle\leq\tfrac{\delta}{2}\|x-y\|^2
\end{align*}
holds for all $x,y\in\reals^n$. This is equivalent to that $\nabla f$ is $\tfrac{1}{\delta}$-cocoercive, see \cite[Theorem~2.1.5]{NesterovLectures} and \cite[Definition~4.4]{bauschkeCVXanal}.
\end{pf}

\subsubsection{Function Properties}
\label{sec:fcn_prop}
\begin{definition}[Strong convexity]
  \label{def:strConv}
  Let $P~:~\reals^n\to\reals^n$ be positive definite. A proper and closed function $f~:~\reals^n\to\overline{\reals}$ is $\sigma$-\emph{strongly convex} w.r.t. $\|\cdot\|_P$ with $\sigma>0$ if $f-\tfrac{\sigma}{2}\|\cdot\|_P^2$ is convex.
\end{definition}
\begin{remark}
If $f$ is differentiable, $\sigma$-strong convexity w.r.t. $\|\cdot\|_P$ can equivalently be defined as that
  \begin{align}
\tfrac{\sigma}{2}\|x-y\|_P^2\leq f(x)-f(y)-\langle \nabla f(y),x-y\rangle
    \label{eq:str_cvx_def}
  \end{align}
holds for all $x,y\in\reals^n$. If $P=\id$, i.e., if the norm is the induced norm, we merely say that $f$ is $\sigma$-strongly convex. If $\sigma=0$, the function is convex.
\end{remark}

There are many smoothness definitions for functions in the literature. We will use the following that implies that the function is in every point majorized and minimized by a norm-squared function.
\begin{definition}[Smoothness]
  Let $P~:~\reals^n\to\reals^n$ be positive semidefinite. A function $f~:~\reals^n\to\reals$ is
  $\beta$-smooth w.r.t. $\|\cdot\|_P$ with $\beta\geq 0$, if it is differentiable and
  \begin{align}
    -\tfrac{\beta}{2}\|x-y\|_P^2\leq f(x)-f(y)-\langle \nabla f(y),x-y\rangle\leq\tfrac{\beta}{2}\|x-y\|_P^2
    \label{eq:fSmooth_cvx}
  \end{align}
  holds for all $x,y\in\reals^n$.
\label{def:smoothness}
\end{definition}

\subsubsection{Connections}
\label{sec:connections}
We will later show that our envelope function satisfies upper and lower bounds of the form
\begin{align}
\tfrac{1}{2}\langle M(x-y),x-y\rangle\leq f(x)-f(y)-\langle\nabla f(y),x-y\rangle\leq\tfrac{1}{2}\langle L(x-y),x-y\rangle
\label{eq:quad_bound_general}
\end{align}
for all $x,y\in\reals^n$ and for different linear operators $M~:~\reals^n\to\reals^n$ and $L~:~\reals^n\to\reals^n$. Depending on $M$ and $L$, we get different properties of $f$ and its gradient $\nabla f$. Some of these are stated below. The results follow immediately from Lemma~\ref{lem:quad_bound_Lipschitz} in Appendix~\ref{app:lemmas} and the definitions of smoothness and strong convexity in Definition~\ref{def:strConv} and Definition~\ref{def:smoothness} respectively.
\begin{proposition}
Assume that $L=-M=\beta I$ with $\beta\geq 0$ in \eqref{eq:quad_bound_general}. Then \eqref{eq:quad_bound_general} is equivalent to that $\nabla f$ is $\beta$-Lipschitz continuous.
\end{proposition}
\begin{proposition}
Assume that $M=\sigma I$ and $L=\beta I$ with $0\leq\sigma\leq \beta$ in \eqref{eq:quad_bound_general}. Then \eqref{eq:quad_bound_general} is equivalent to that $\nabla f$ is $\beta$-Lipschitz continuous and $f$ is $\sigma$-strongly convex.
\end{proposition}
\begin{proposition}
Assume that $L=-M$ and that $L$ is positive definite. Then \eqref{eq:quad_bound_general} is equivalent to that $f$ is $1$-smooth w.r.t. $\|\cdot\|_L$.
\end{proposition}
\begin{proposition}
Assume that $M$ and $L$ are positive definite. Then \eqref{eq:quad_bound_general} is equivalent to that $f$ is $1$-smooth w.r.t. $\|\cdot\|_L$ and $1$-strongly convex w.r.t. $\|\cdot\|_M$.
\end{proposition}

\section{Envelope Functions}

To find a fixed-point of a nonexpansive mapping $S$ using an averaged iteration of that mapping, is the basis for many first-order optimization methods. Based on ideas from \cite{Panos_fb_newton,Panos_acc_DR_2014}, we present another method to find such a fixed-point. We create an envelope function whose stationary points coincide with the fixed-points of the operator $S$. For forward-backward splitting and Douglas-Rachford splitting, such envelopes have been proposed in \cite{Panos_fb_newton} and \cite{Panos_acc_DR_2014} respectively. These envelope functions turn out to be special cases of the envelopes we propose, see Section~\ref{sec:special_cases}. The envelope functions often possess favorable properties such as convexity and Lipschitz continuity of the gradient. Then, any method to find a stationary point (in the convex case, a minimizer) of the envelope function can be used to find a fixed-point to the nonexpansive mapping $S$.

To formulate our envelope function, we assume that the nonexpansive operator $S$ is a composition of $S_2$ and $S_1$, i.e., $S=S_2S_1$. We make the following basic assumptions on $S_1$ and $S_2$, that sometimes will be sharpened or relaxed:
\begin{assumption}
Suppose that:
\begin{enumerate}[(i)]
\item $S_1~:~\reals^n\to\reals^n$ and $S_2~:~\reals^n\to\reals^n$ are nonexpansive
\item $S_1=\nabla f_1$ and $S_2=\nabla f_2$ for some differentiable functions $f_1~:~\reals^n\to\reals$ and $f_2~:~\reals^n\to\reals$
\item $S_1~:~\reals^n\to\reals^n$ is affine, i.e., $S_1x=Px+q$ and $f_1(x)=\tfrac{1}{2}\langle Px,x\rangle+\langle q,x\rangle$, where $P\in\reals^{n\times n}$ is a a self-adjoint nonexpansive linear operator and $q\in\reals^n$
\end{enumerate}
\label{ass:Si_prop}
\end{assumption}
\begin{remark}
Part (iii) of the assumption means that $P$ is symmetric with eigenvalues in the interval $[-1,1]$.
\end{remark}
Now, we are ready to define the general envelope function whose properties we will investigate in this paper:
\begin{align}
F(x):=\tfrac{1}{2}\langle Px,x\rangle- f_2(\nabla f_1(x)).
\label{eq:env}
\end{align}
The gradient of this function is given by
\begin{align}
\nabla F(x) = Px-\nabla^2f_1(x)\nabla f_2(\nabla f_1(x))=Px-PS_2(S_1x)=P(x-S_2S_1x).
\label{eq:env_grad}
\end{align}
The set of stationary points to the envelope function $F$ is the set of points for which the gradient is zero. This set is denoted as follows:
\begin{align}
X^\star:=\{x~|~\nabla F(x)=0\}.
\label{eq:env_stat_points}
\end{align}


\subsection{Basic Properties of the Envelope Function}

\label{sec:basic_prop}

Here, we list some basic properties of the envelope function \eqref{eq:env}. The first two results are special cases and direct corollaries of a more general result in Theorem~\ref{thm:env_finer_prop}, and therefore not proven here.
\begin{proposition}
Suppose that Assumption~\ref{ass:Si_prop} holds. Then the gradient of $F$ is 2-Lipschitz continuous. That is, $\nabla F$ satisfies
\begin{align*}
\|\nabla F(x)-\nabla F(y)\|\leq 2\|x-y\|
\end{align*}
for all $x,y\in\reals^n$.
\label{prp:env_smooth}
\end{proposition}
\begin{proposition}
Suppose that Assumption~\ref{ass:Si_prop} holds and that $P$, the operator defining the linear part of $S_1$, is positive semidefinite. Then $F$ is convex.
\label{prp:env_conv}
\end{proposition}
So, if $P$ is positive semidefinite, then the envelope function $F$ is convex and differentiable with a Lipschitz continuous gradient. The set of stationary points of $F$ also has a close relationship with the fixed-points of $S=S_2S_1$. This is shown next.
\begin{proposition}
Suppose that Assumption~\ref{ass:Si_prop} holds and that $P$ is nonsingular. Then $X^\star={\rm{fix}}(S_2S_1)$ where $X^\star$ is defined in \eqref{eq:env_stat_points} and the fixed-point set ${\rm{fix}}(S_2S_1)$ is ${\rm{fix}}(S_2S_1)=\{x\in\reals^n:S_2S_1x=x\}$. If in addition $P$ is positive definite, then $\Argmin_x F(x)=X^\star={\rm{fix}}(S_2S_1)$.
\label{prp:fp_sp_agree}
\end{proposition}
\begin{pf}
The first claim follows directly from \eqref{eq:env_grad}. The second claim follows from \eqref{eq:env_grad} and that $F$ is convex when $P$ is positive (semi)definite, see Proposition~\ref{prp:env_conv}.
\end{pf}
These three results show that if $P$ is positive definite, a fixed-point to $S_2S_1$ can be found by minimizing the differentiable convex function $F$, which has a 2-Lipschitz continuous gradient.


\subsection{Finer Properties of the Envelope Function}

\label{sec:finer_prop}

Here, we establish some finer properties of the envelope function. We start with a general result on upper and lower bounds for the envelope function. This result uses stronger assumptions on $S_2$ than nonexpansiveness, namely that it is $\alpha$-averaged and $\beta$-negatively averaged with $\alpha,\beta\in(0,1]$, see Definition~\ref{def:avg} and Definition~\ref{def:negavg}. We state this as an assumption.
\begin{assumption}
The operator $S_2$ is $\alpha$-averaged and $\beta$-negatively averaged with $\alpha\in(0,1]$ and $\beta\in(0,1]$.
\label{ass:S2_finer_prop}
\end{assumption}

\begin{theorem}
Suppose that Assumption~\ref{ass:Si_prop} and Assumption~\ref{ass:S2_finer_prop} hold. Further, let $\delta_{\alpha}=2\alpha-1$ and $\delta_{\beta}=2\beta-1$. Then the envelope function $F$ in \eqref{eq:env} satisfies
\begin{align*}
 F(x)-F(y)-\langle\nabla F(y),x-y\rangle\geq\tfrac{1}{2} \langle (P-\delta_{\beta}P^2)(x-y),x-y\rangle
\end{align*}
and
\begin{align*}
F(x)-F(y)-\langle\nabla F(y),x-y\rangle\leq \tfrac{1}{2}\langle(P+\delta_{\alpha}P^2)(x-y),x-y\rangle
\end{align*}
for all $x,y\in\reals^n$.
\label{thm:env_finer_prop}
\end{theorem}
A proof to this result is found in Appendix~\ref{app:env_finer_prop_pf}.

As seen in Section~\ref{sec:connections}, such bounds have many implications on the properties of the function. Next, we provide some in the form of corollaries.
\begin{corollary}
Suppose that Assumption~\ref{ass:Si_prop} and Assumption~\ref{ass:S2_finer_prop} hold and that $P$ is positive semidefinite. Let $\delta_{\alpha}=2\alpha-1$ and $\delta_{\beta}=2\beta-1$. Then
\begin{align*}
\tfrac{1}{2}\|x-y\|_{P-\delta_{\beta}P^2}^2\leq F(x)-F(y)-\langle\nabla F(y),x-y\rangle\leq \tfrac{1}{2}\|x-y\|_{P+\delta_{\alpha}P^2}^2
\end{align*}
where $P-\delta_{\beta} P^2$ is positive semidefinite.
\label{cor:quad_lower_upper_scaled}
\end{corollary}
\begin{pf}
It follows directly from Theorem~\ref{thm:env_finer_prop} and Lemma~\ref{lem:smallest_eigenvalue} in Appendix~\ref{app:lemmas}.
\end{pf}

\begin{corollary}
Suppose that Assumption~\ref{ass:Si_prop} and Assumption~\ref{ass:S2_finer_prop} hold and that either of the following holds:
\begin{enumerate}[(i)]
\item $P$ is positive definite and contractive
\item $P$ is positive definite and $\beta\in(0,1)$ in the negative averagedness
\end{enumerate}
Let $\delta_{\alpha}=2\alpha-1$ and $\delta_{\beta}=2\beta-1$. Then $F$ is 1-strongly convex w.r.t. $\|\cdot\|_{P-\delta_{\beta}P^2}$ and 1-smooth w.r.t. $\|\cdot\|_{P+\delta_{\alpha} P^2}$.
\label{cor:str_conv_smooth_general}
\end{corollary}
\begin{pf}
To show the strong convexity claim, it is sufficient to apply Theorem~\ref{thm:env_finer_prop} and show that $P-\delta_{\beta}P^2$ is positive definite, i.e., that $\lambda_{\min}(P-\delta_{\beta}P^2)$ is positive. In {\it{(i)}}, $\lambda_i(P)\in(0,1)$ and $\delta_{\beta}\in(-1,1]$ and in {\it{(ii)}}, $\lambda_i(P)\in(0,1]$ and $\delta_{\beta}\in(-1,1)$. From Lemma~\ref{lem:smallest_eigenvalue} it follows that in both cases, $\lambda_{\min}(P-\delta_{\beta}P^2)$ is positive. The smoothness claim follows immediately from Theorem~\ref{thm:env_finer_prop} and Definition~\ref{def:smoothness}.
\end{pf}
Next, we show a less tight characterization of the envelope function that does not take the shape of the upper and lower bounds into account.
\begin{corollary}
Suppose that Assumption~\ref{ass:Si_prop} and Assumption~\ref{ass:S2_finer_prop} hold. Let $m=\lambda_{\min}(P)$, $L=\lambda_{\max}(P)$, $\delta_{\alpha}=2\alpha-1\in[-0.5,1]$, and $\delta_{\beta}=2\beta-1\in[-0.5,1]$. Then
\begin{align*}
\tfrac{\beta_l}{2}\|x-y\|^2\leq F(x)-F(y)-\langle\nabla F(y),x-y\rangle\leq \tfrac{\beta_u}{2}\|x-y\|^2
\end{align*}
where $\beta_l=\min(m(1-\delta_{\beta}m),L(1-\delta_{\beta}L))$ and $\beta_u=L(1+\delta_{\alpha}L)$.
\label{cor:quad_upper_lower_id}
\end{corollary}
\begin{pf}
This follows from Theorem~\ref{thm:env_finer_prop}, Lemma~\ref{lem:smallest_eigenvalue}, and Lemma~\ref{lem:largest_eigenvalue}.
\end{pf}
We restricted $\delta_{\alpha}$ and $\delta_{\beta}$ to $[-0.5,1]$ (i.e, $\alpha$ and $\beta$ to $[0.25,1]$) in this result for convenience of the statement. Similar results for other $\delta_{\beta}$ and $\delta_{\alpha}$ (and a sharpening of the result when $\delta_{\beta}\in[-0.5,0]$) can be concluded from Lemma~\ref{lem:smallest_eigenvalue} and Lemma~\ref{lem:largest_eigenvalue}.

From Corollary~\ref{cor:quad_upper_lower_id}, the following two results are immediate.
\begin{corollary}
Suppose that Assumption~\ref{ass:Si_prop} and Assumption~\ref{ass:S2_finer_prop} hold. Let $\delta_{\alpha}=2\alpha-1\in[-0.5,1]$, $\delta_{\beta}=2\beta-1\in[-0.5,1]$, $m=\lambda_{\min}(P)$, and $L=\lambda_{\max}(P)$ and suppose that either of the following two conditions holds: 
\begin{enumerate}[(i)]
\item $P$ is positive definite with $\lambda_{\min}(P)\in(0,1)$ and $\lambda_{\max}(P)\in[m,1)$
\item $P$ is positive definite with $\lambda_{\min}(P)\in(0,1]$ and $\delta_{\beta}=2\beta-1\in[-0.5,1)$
\end{enumerate}
Then $F$ is $\min(m(1-\delta_{\beta}m),L(1-\delta_{\beta}L))$-strongly convex (w.r.t. $\|\cdot\|$) and $L(1+\delta_{\alpha}L)$-smooth (w.r.t. $\|\cdot\|$).
\label{cor:str_conv_smooth_id}
\end{corollary}
\begin{corollary}
Suppose that Assumption~\ref{ass:Si_prop} and Assumption~\ref{ass:S2_finer_prop} hold and that $P$ is positive semidefinite, i.e., that $\lambda_{\min}(P)\geq 0$. Let $L=\lambda_{\max}(P)$, $\delta_{\beta}=2\beta-1\in[-0.5,1]$, and $\delta_{\alpha}=2\alpha-1\in[-0.5,1]$. Then $F$ is convex and it is $L(1+\delta_{\alpha}L)$-smooth (or equivalently $\nabla F$ is $L(1+\delta_{\alpha} L)$-Lipschitz continuous).
\label{cor:env_conv_smooth}
\end{corollary}

The results in Theorem~\ref{thm:env_finer_prop} and its corollaries hold for $\alpha$-averaged and $\beta$-negatively averaged operators $S_2$. In Proposition~\ref{prp:prop_relation}, some properties that are equivalent to averagedness and negative averagedness are stated. Therefore, we can use these equivalent properties instead when stating the above results. This is done in the following to propositions.
\begin{proposition}
Suppose that Assumption~\ref{ass:Si_prop} holds and that $S_2$ is $\delta$-Lipschitz continuous with $\delta\in[0,1]$. Then all results in this section hold with $\delta_{\beta}=\delta_{\alpha}=\delta$.
\label{prp:delta_when_lipschitz}
\end{proposition}
\begin{proposition}
Suppose that Assumption~\ref{ass:Si_prop} holds and that $S_2$ is $\tfrac{1}{\delta}$-cocoercive with $\delta\in[0,1]$. Then all results in this section hold with $\delta_{\beta}=\delta$ and $\delta_{\alpha}=0$.
\label{prp:delta_when_coco}
\end{proposition}

\subsection{Relation to Averaged Operator Iteration}

\label{sec:rel_to_avg}

As noted in \cite{Panos_fb_newton,Panos_acc_DR_2014}, the forward-backward and Douglas-Rachford splitting methods are variable metric gradient methods applied to their respective envelope functions. In our setting with $S_1$ being affine, it reduces to a fixed-metric scaled gradient method. Here, we show that this observation holds also in our setting.

We apply the following scaled gradient method to the envelop function $F$:
\begin{align*}
x^{k+1}&=x^k-\alpha P^{-1}\nabla F(x^k).
\end{align*}
This gives
\begin{align*}
x^{k+1}&=x^k-\alpha P^{-1}\nabla F(x^k)\\
&=x^k-\alpha P^{-1}P(S_2S_1x^k-x^k)\\
&=x^k-\alpha(S_2S_1x^k-x^k)\\
&=(1-\alpha)x^k+\alpha S_2S_1x^k,
\end{align*}
which is an averaged iteration of the nonexpansive mapping $S_2S_1$ for $\alpha\in(0,1)$. Therefore, the basic averaged iteration can be interpreted as a scaled gradient method applied to the envelope function.

This is most probably not the most efficient way to find a stationary point of the envelope function (or equivalently a fixed-point to $S_2S_1$). At least in the convex setting (for the envelope), there are numerous alternative methods that can minimize smooth functions such as truncated Newton methods, quasi-Newton methods, and nonlinear conjugate gradient descent. See \cite{Nocedal} for an overview of such methods and \cite{Panos_fb_newton,Panos_fb_quasi_newton} for some of these methods applied to the forward-backward envelope. Evaluating which ones that are most efficient and devising new methods to improve performance is outside the scope of this paper.

\section{Special Cases}

\label{sec:special_cases}

In this section, we present a generalization of the envelope function in the previous section. This envelope has four known special cases, namely the Moreau envelope~\cite{Moreau1965}, the forward-backward envelope \cite{Panos_fb_newton,Panos_fb_quasi_newton}, the Douglas-Rachford envelope \cite{Panos_acc_DR_2014}, and the ADMM envelope (which is a special case of the Douglas-Rachford envelope). 

The generalization incorporates envelopes for iterations where $f_1$ that defines $S_1$ through $S_1=\nabla f_1$ is twice continuously differentiable (as opposed to quadratic in the previous section). The more general envelope function is
\begin{align}
  F(x)=\langle\nabla f_1(x),x\rangle-f_1(x)-f_2(\nabla f_1(x)).
\label{eq:env_general}
\end{align}
When $f_1(x)=\tfrac{1}{2}\langle Px,x\rangle+\langle q,x\rangle$ it reduces to \eqref{eq:env} since then 
\begin{align*}
\langle\nabla f_1(x),x\rangle-f_1(x)=\langle Px+q,x\rangle-(\tfrac{1}{2}\langle Px,x\rangle+\langle q,x\rangle)=\tfrac{1}{2}\langle Px,x\rangle.
\end{align*}
The gradient of the envelope function in \eqref{eq:env_general} is
\begin{align*}
\nabla F(x)&=\nabla^2f_1(x)x+\nabla f_1(x)-\nabla f_1(x)-\nabla^2f_1(x)\nabla f_2(\nabla f_1(x))\\
&=\nabla^2f_1(x)(x-\nabla f_2(\nabla f_1(x)))\\
&=\nabla^2f_1(x)(x-S_2S_1x).
\end{align*}
If $\nabla^2f_1(x)$ is nonsingular for all $x$, the set of stationary points of the envelope coincides with the fixed-point set of $S=S_2S_1$. We do not provide any properties of the envelope functions in this setting (it is left as future work), but merely show that that it generalizes the previously known envelope functions.

In the more restricted setting with $S_1=\nabla f_1$ being affine, we provide envelope function properties that coincide with or sharpen corresponding results in the literature for the special cases.

\subsection{Preliminaries}

Before we present the special cases, we introduce some functions whose gradients are operators that are used in the respective underlying methods. Most importantly, we will introduce a function whose gradient is the proximal operator, which is defined as follows:
\begin{align*}
  \prox_{\gamma f}(z):=\argmin_{x}\{f(x)+\tfrac{1}{2\gamma}\|x-z\|^2\},
\end{align*}
where $\gamma>0$ is a parameter. To do this, we introduce the following function which is a scaling and regularization of $f$:
\begin{align}
  r_{\gamma f}(x):=\gamma f(x)+\tfrac{1}{2}\|x\|^2
\label{eq:fcn_grad_prox}
\end{align}
This is related to the proximal operator of $f$ as follows:
\begin{proposition}
Suppose that $f~:~\reals^n\to\reals\cup\{\infty\}$ is proper closed and convex and that $\gamma>0$. The proximal operator ${\rm{prox}}_{\gamma f}$ then satisfies
\begin{align*}
{\rm{prox}}_{\gamma f} = \nabla r_{\gamma f}^*
\end{align*}
where $r_{\gamma f}$ is defined in \eqref{eq:fcn_grad_prox}.
\end{proposition}
This result is from \cite[Theorem~31.5,~Theorem~16.4]{Rockafellar} and implies that the proximal operator is the gradient of a convex function.

A special case is when $f=\iota_C$, where $\iota_C$ is the indicator function for the nonempty closed and convex set $C$. The proximal operator then reduces to the projection operator. The projection operator onto $C$ is denoted by $\Pi_C$ and the corresponding regularized function is denoted and defined by
\begin{align}
  r_{C}(x):=\iota_{C}(x)+\tfrac{1}{2}\|x\|^2.
\label{eq:fcn_grad_proj}
\end{align}
With this notation, $\Pi_{C}(x)=\nabla r_{C}^*(x)$. Next, we introduce a linear combination between $r^*$ and $\tfrac{1}{2}\|\cdot\|^2$, namely
\begin{align}
p_{\gamma f}^{\alpha}(x):=\alpha r_{\gamma f}^*(x)+\tfrac{1-\alpha}{2}\|x\|^2,
\label{eq:fcn_grad_rel_prox}
\end{align}
where we typically require that $\alpha\in(0,2]$. The gradient of $p_{\gamma f}^{\alpha}$ is denoted by $P_{\gamma f}^{\alpha}$ and is given by
\begin{align}
P_{\gamma f}^{\alpha}(x):=\nabla p_{\gamma f}^{\alpha}(x)=\alpha{\rm{prox}}_{\gamma f}(x)+(1-\alpha)x.
\label{eq:rel_prox}
\end{align}
This is called a relaxed proximal mapping. Some special cases of this will have their own notation. Letting $\alpha=2$, we get the reflected proximal operator 
\begin{align}
R_{\gamma f}(x):=P_{\gamma f}^2(x)=2{\rm{prox}}_{\gamma f}(x)-x.
\label{eq:refl_prox}
\end{align}
When $f=\iota_C$, we will use notation $p_{C}^{\alpha}$, $P_{C}^{\alpha}$, and $R_C$ for \eqref{eq:fcn_grad_rel_prox}, \eqref{eq:rel_prox}, and \eqref{eq:refl_prox} respectively. That is
\begin{align}
\label{eq:fcn_grad_rel_proj} p_{C}^{\alpha}(x)&:=\alpha r_{C}^*(x)+\tfrac{1-\alpha}{2}\|x\|^2,\\
\label{eq:rel_proj} P_{C}^{\alpha}(x)&:=\nabla p_{C}^{\alpha}(x)=\alpha\Pi_{C}(x)+(1-\alpha)x\\
\label{eq:refl_proj} R_C(x)&:=2\Pi_C(x)-x.
\end{align}
We refer to \eqref{eq:rel_proj} as a relaxed projection, and \eqref{eq:refl_proj} as a reflection. So, the proximal and projected operators and their relaxed and reflected variants are gradients of functions.

We conclude with the straightforward observation that
\begin{align*}
(x-\gamma\nabla f(x))=\nabla\left(\tfrac{1}{2}\|x\|^2-\gamma f(x)\right).
\end{align*}
That is, the gradient step operator is the gradient of the function $\tfrac{1}{2}\|x\|^2-\gamma f(x)$.

\subsection{The Proximal Point Algorithm}

The proximal point algorithm solves problems of the form
\begin{align*}
{\hbox{minimize }} f(x)
\end{align*}
where $f~:~\reals^n\to\reals\cup\{\infty\}$ is proper closed and convex.

The algorithm repeatedly applies the proximal operator of $f$ and is given by
\begin{align}
x^{k+1} = \prox_{\gamma f}(x^k),
\label{eq:ppa}
\end{align}
where $\gamma>0$ is a parameter. This algorithm is mostly of conceptual interest since it is often as computationally demanding to evaluate the prox as to minimize the function $f$ itself.

Its envelope function, which is called the Moreau envelope \cite{Moreau1965}, is a scaled version of our envelope $F$ in \eqref{eq:env}. The scaling factor is $\gamma^{-1}$ and $F$ in \eqref{eq:env} is obtained by letting $S_1x=\nabla f_1(x)=x$, i.e., $P=\id$ and $q=0$, and $f_2=r_{\gamma f}^*$, where $r_{\gamma f}$ is defined in \eqref{eq:fcn_grad_prox}. The resulting envelope function $f^{\gamma}$ is given by
\begin{align}
f^{\gamma}(x)=\gamma^{-1}F(x)=\gamma^{-1}\left(\tfrac{1}{2}\|x\|^2-r_{\gamma f}^{*}(x)\right),
\label{eq:Moreau_env}
\end{align}
and its gradient satisfies
\begin{align*}
\nabla f^{\gamma}(x)=\gamma^{-1}\left(x-\prox_{\gamma f}(x)\right).
\end{align*}
The following properties of the Moreau envelope follow directly from Corollary~\ref{cor:env_conv_smooth} and Proposition~\ref{prp:delta_when_coco} since the proximal operator is 1-cocoercive (see Remark~\ref{rem:coco_avg} and \cite[Proposition~12.27]{bauschkeCVXanal}). 
\begin{proposition}
The Moreau envelope $f^{\gamma}$ in \eqref{eq:Moreau_env} is differentiable and convex and $\nabla f^{\gamma}$ is $\gamma^{-1}$-Lipschitz continuous.
\end{proposition}
This coincides with previously known properties of the Moreau envelope, see \cite[Chapter~12]{bauschkeCVXanal}.

\subsection{Forward-Backward Splitting}

Forward-backward splitting solves problems of the form
\begin{align}
{\hbox{minimize }} f(x)+g(x)
\label{eq:FB_prob}
\end{align}
where $f~:~\reals^n\to\reals$ is convex with an $L$-Lipschitz (or equivalently $\tfrac{1}{L}$-cocoercive) gradient, and $g~:~\reals^n\to\reals\cup\{\infty\}$ is proper closed and convex.

The algorithm performs a forward step then a backward step and is given by
\begin{align}
x^{k+1}=\prox_{\gamma g}(\id-\gamma\nabla f)x^k,
\label{eq:FB}
\end{align}
where $\gamma\in(0,\tfrac{2}{L})$ is a parameter. 

The envelope function, which is called the forward-backward envelope \cite{Panos_fb_newton,Panos_fb_quasi_newton}, is a scaled version of our envelope $F$ in \eqref{eq:env_general} and applies when $f$ is twice continuously differentiable and $\nabla F$ is Lipschitz continuous. The scaling factor is $\gamma^{-1}$ and $F$ in \eqref{eq:env_general} is obtained by letting $f_1=\tfrac{1}{2}\|\cdot\|^2-\gamma f$ and $f_2=r_{\gamma g}^*$, where $r_{\gamma g}$ is defined in \eqref{eq:fcn_grad_prox}. The resulting forward-backward envelope function is
\begin{align*}
F_{\gamma}^{\rm{FB}}(x)=\gamma^{-1}\left(\langle x-\gamma\nabla f(x),x\rangle-(\tfrac{1}{2}\|x\|^2-\gamma f(x))-r_{\gamma g}^*(x-\gamma\nabla f(x))\right).
\end{align*}
The gradient of this function is
\begin{align*}
\nabla F_{\gamma}^{\rm{FB}}(x)&=\gamma^{-1}\big((\id-\gamma\nabla^2 f(x))x+(x-\gamma\nabla f(x))-(x-\gamma\nabla f(x))\\
&\quad-(\id-\gamma\nabla^2 f(x))\prox_{\gamma g}(x-\gamma\nabla f(x))\big)\\
&=\gamma^{-1}(\id-\gamma\nabla^2 f(x))\left(x-\prox_{\gamma g}(x-\gamma\nabla f(x))\right)
\end{align*}
which coincides with the gradient in \cite{Panos_fb_newton,Panos_fb_quasi_newton}. As described in \cite{Panos_fb_newton,Panos_fb_quasi_newton}, the stationary points of the envelope coincide with the fixed-points of $x-\prox_{\gamma g}(x-\gamma\nabla f(x))$ if $(\id-\gamma\nabla^2 f(x))$ is nonsingular.

\subsubsection{$S_1$ affine}

We provide properties of the forward-backward envelope in the more restrictive setting where $S_1=\nabla f_1=(\id-\gamma \nabla f)$ is affine. This happens if $f$ is convex quadratic, i.e., $f(x)=\tfrac{1}{2}\langle Hx,x\rangle+\langle h,x\rangle$ with $H\in\reals^{n\times n}$ positive semidefinite  and $h\in\reals^n$. Then $S_1x=Px+q$ with $P=(\id-\gamma H)$ and $q=-\gamma h$.

In this setting, the following result follows immediately from Corollary~\ref{cor:quad_lower_upper_scaled} and Proposition~\ref{prp:delta_when_coco} (where Proposition~\ref{prp:delta_when_coco} is invoked since $S_2=\prox_{\gamma g}$ is 1-cocoercive, see Remark~\ref{rem:coco_avg} and \cite[Proposition~12.27]{bauschkeCVXanal}).
\begin{proposition}
Assume that $f(x)=\tfrac{1}{2}\langle Hx,x\rangle+\langle h,x\rangle$ and $\gamma\in(0,\tfrac{1}{L})$ where $L=\lambda_{\max}(H)$. Then the forward-backward envelope $F_{\gamma}^{\rm{FB}}$ satisfies
\begin{align*}
\tfrac{1}{2\gamma}\|x-y\|_{P-P^2}^2&\leq F_{\gamma}^{\rm{FB}}(x)-F_{\gamma}^{\rm{FB}}(y)-\langle\nabla F_{\gamma}^{\rm{FB}}(y),x-y\rangle\leq \tfrac{1}{2\gamma}\|x-y\|_P^2
\end{align*}
for all $x,y\in\reals^n$, where $P=(\id-\gamma H)$ is positive definite. If in addition $\lambda_{\min}(H)=m>0$, then $P-P^2$ is positive definite and $F_{\gamma}^{\rm{FB}}$ is $\gamma^{-1}$-strongly convex w.r.t. $\|\cdot\|_{P-P^2}$.
\label{prp:fb_sharp}
\end{proposition}
Less tight bounds for the forward-backward envelope are provided next. These follow immediately from Corollary~\ref{cor:str_conv_smooth_id}, Corollary~\ref{cor:env_conv_smooth}, and Proposition~\ref{prp:delta_when_coco}.
\begin{proposition}
Assume that $f(x)=\tfrac{1}{2}\langle Hx,x\rangle+\langle h,x\rangle$, that $\gamma\in(0,\tfrac{1}{L})$ where $L=\lambda_{\max}(H)$, and that $m=\lambda_{\min}(H)\geq 0$. Then the forward-backward envelope $F_{\gamma}^{\rm{FB}}$ is $\gamma^{-1}(1-\gamma m)$-smooth and $\min\left((1-\gamma m)m,(1-\gamma\L)L\right)$-strongly convex (both w.r.t. to the induced norm $\|\cdot\|$).
\end{proposition}
This result is a less tight version of Proposition~\ref{prp:fb_sharp}, but is a slight improvement of the corresponding result in \cite[Theorem~2.3]{Panos_fb_newton}. The strong convexity moduli are the same, but this smoothness constant is a factor two smaller.



\subsection{Douglas-Rachford Splitting}

\label{sec:DR}
Douglas-Rachford splitting solves problems of the form
\begin{align}
{\hbox{minimize }} f(x)+g(x)
\label{eq:DR_prob}
\end{align}
where $f~:~\reals^n\to\reals\cup\{\infty\}$ and $g~:~\reals^n\to\reals\cup\{\infty\}$ are proper closed and convex functions.

The algorithm performs two reflection steps \eqref{eq:refl_prox}, then an averaging according to
\begin{align}
z^{k+1}=(1-\alpha)z^k+\alpha R_{\gamma g}R_{\gamma f}z^k
\label{eq:DR}
\end{align}
where $\gamma>0$ and $\alpha\in(0,1)$ are parameters. The objective is to find a fixed-point $\bar{z}$ to $R_{\gamma g}R_{\gamma f}$, from which a solution to \eqref{eq:DR_prob} can be computed as $\prox_{\gamma f}\bar{z}$, see \cite[Proposition~25.1]{bauschkeCVXanal}.

The envelope function from \cite{Panos_acc_DR_2014}, which is called the Douglas-Rachford envelope, is a scaled version of the basic envelope function $F$ in \eqref{eq:env_general} and applies when $f$ is twice continuously differentiable and $\nabla F$ is Lipschitz continuous. The scaling factor is $(2\gamma)^{-1}$ and $F$ is obtained by letting $f_1=p_{\gamma f}^2$ with gradient $\nabla f_1=S_1=R_{\gamma f}$ and $f_2 = p_{\gamma g}^2$, where $p_{\gamma g}^2$ is defined in \eqref{eq:fcn_grad_rel_prox}. The Douglas-Rachford envelope function becomes
\begin{align}
F_{\gamma}^{\rm{DR}}(z)=(2\gamma)^{-1}\left(\langle R_{\gamma f}(z),z\rangle-p_{\gamma f}^2(z)-p_{\gamma g}^2(R_{\gamma f}z)\right).
\label{eq:DR_env}
\end{align}
The gradient of this function is
\begin{align*}
\nabla F_{\gamma}^{\rm{DR}}(z)&=(2\gamma)^{-1}\big(\nabla R_{\gamma f}(z)z+R_{\gamma f}-R_{\gamma f}-\nabla R_{\gamma f}(z)R_{\gamma g}(R_{\gamma f}(z))\big)\\
&=(2\gamma)^{-1}\nabla R_{\gamma f}(z)(z-R_{\gamma g}R_{\gamma f}(z)).
\end{align*}
which coincides with the gradient in \cite{Panos_acc_DR_2014} since $\nabla R_{\gamma f}=2\nabla \prox_{\gamma f}-\id$ and 
\begin{align*}
z-R_{\gamma g}R_{\gamma f}z&=z-2\prox_{\gamma g}(2\prox_{\gamma f}(z)-z)+2\prox_{\gamma f}(z)-z\\
&=2(\prox_{\gamma f}(z)-\prox_{\gamma g}(2\prox_{\gamma f}(z)-z)).
\end{align*}
As described in \cite{Panos_acc_DR_2014}, the stationary points of the envelope coincide with the fixed-points of $x-R_{\gamma g}R_{\gamma f}$ if $\nabla R_{\gamma f}$ is nonsingular.

\subsubsection{$S_1$ affine}

We state properties of the Douglas-Rachford envelope in the more restrictive setting where $S_1=R_{\gamma f}$ is affine. This holds if $f$ is convex quadratic, i.e., of the form
\begin{align*}
f(x)=\tfrac{1}{2}\langle Hx,x\rangle+\langle h,x\rangle.
\end{align*}
The operator $S_1$ becomes
\begin{align*}
S_1(z) = R_{\gamma f}(z) = 2(\id+\gamma H)^{-1}(z-\gamma h)-z,
\end{align*}
which confirms that it is affine. We implicitly define $P$ and $q$ through $S_1=R_{\gamma f}=P(\cdot)+q$, and note that they are given by $P=2(\id+\gamma H)^{-1}-\id$ and $q=-2\gamma (\id+\gamma H)^{-1} h$.

In this setting, the following result follows immediately from Corollary~\ref{cor:quad_lower_upper_scaled} since $S_2=R_{\gamma g}$ is nonexpansive (1-averaged and 1-negatively averaged).
\begin{proposition}
Assume that $f(x)=\tfrac{1}{2}\langle Hx,x\rangle+\langle h,x\rangle$ and $\gamma\in(0,\tfrac{1}{L})$ where $L=\lambda_{\max}(H)$. Then the Douglas-Rachford envelope $F_{\gamma}^{\rm{DR}}$ satisfies
\begin{align*}
\tfrac{1}{4\gamma}\|z-y\|_{P-P^2}^2&\leq F_{\gamma}^{\rm{DR}}(z)-F_{\gamma}^{\rm{DR}}(z)-\langle\nabla F_{\gamma}^{\rm{DR}}(y),z-y\rangle\leq \tfrac{1}{4\gamma}\|z-y\|_{P+P^2}^2
\end{align*}
for all $y,z\in\reals^n$, where $P=2(\id+\gamma H)^{-1}-\id$ is positive definite. If in addition $\lambda_{\min}(H)=m>0$, then $P-P^2$ is positive definite and $F_{\gamma}^{\rm{DR}}$ is $(2\gamma)^{-1}$-strongly convex w.r.t. $\|\cdot\|_{P-P^2}$.
\label{prp:DR_bounds_general}
\end{proposition}

The following less tight characterization of the Douglas-Rachford envelope follows from Corollary~\ref{cor:str_conv_smooth_id} and Corollary~\ref{cor:env_conv_smooth}.
\begin{proposition}
Assume that $f(x)=\tfrac{1}{2}\langle Hx,x\rangle+\langle h,x\rangle$, that $\gamma\in(0,\tfrac{1}{L})$ where $L=\lambda_{\max}(H)$, and that $m=\lambda_{\min}(H)\geq 0$. Then the Douglas-Rachford envelope $F_{\gamma}^{\rm{DR}}$ is $\tfrac{1-\gamma m}{(1+\gamma m)^2}\gamma^{-1}$-smooth and $\min\left(\tfrac{(1-\gamma m) m}{(1+\gamma m)^2},\tfrac{(1-\gamma L)L}{(1+\gamma L)^2}\right)$-strongly convex.
\label{prp:DR_bounds_id}
\end{proposition}

This result is more conservative than the one in Proposition~\ref{prp:DR_bounds_general}, but improves on \cite[Theorem~2]{Panos_acc_DR_2014}. The strong convexity modulus coincides with the corresponding one in \cite[Theorem~2]{Panos_acc_DR_2014}. The smoothness constant is $\tfrac{1}{1+\gamma m}$ times that in \cite[Theorem~2]{Panos_acc_DR_2014}, i.e., it is slightly smaller.



\subsection{ADMM}

The alternating direction method of multipliers (ADMM) solves problems of the form \eqref{eq:DR_prob}. It is well known \cite{Gabay83} that ADMM can be interpreted as Douglas-Rachford applied to the dual of \eqref{eq:DR_prob}, namely to
\begin{align}
{\hbox{minimize }} f^*(\mu)+g^*(-\mu).
\label{eq:ADMM_prob}
\end{align}
So the algorithm is given by
\begin{align}
v^{k+1} = (1-\alpha)v^k+\alpha R_{\rho (g^*\circ-\id)}R_{\rho f}v^k
\label{eq:ADMM}
\end{align}
where $\rho>0$ is a parameter, and $R_{\rho f}$ the reflected proximal operator \eqref{eq:refl_prox} and $(g^*\circ-\id)$ is the composition that satisfies $(g^*\circ-\id)(\mu)=g^*(-\mu)$.

In accordance with the Douglas-Rachford envelope \eqref{eq:DR_env}, the ADMM envelope is defined as
\begin{align}
F_{\rho}^{\rm{ADMM}}(v)=(2\rho)^{-1}\left(\langle R_{\rho f^*}(v),v\rangle-p_{\rho f^*}^2(v)-p_{\rho (g^*\circ -\id)}^2(R_{\rho f^*}v)\right).
\label{eq:DR_ADMM}
\end{align}
and its gradient becomes
\begin{align*}
\nabla F_{\rho}^{\rm{ADMM}}(v)=(2\rho)^{-1}\nabla R_{\rho f^*}(v)(v-R_{\rho (g^*\circ -\id)}R_{\rho f^*}(v)).
\end{align*}

In this section, we relate the ADMM algorithm and its envelope function to the Douglas-Rachford counterparts. To do so, we need the following lemma which is proven in Appendix~\ref{app:refl_prox_conj_relations}.
\begin{lemma}
Let $g~:~\reals^n\to\reals\cup\{\infty\}$ and be proper closed and convex and $\rho>0$. Then
\begin{align*}
R_{\rho g^*}(x) &= -\rho R_{\rho^{-1}g}(\rho^{-1}x)\\
R_{\rho (g^*\circ -\id)}(x) &= \rho R_{\rho^{-1} g}(-\rho^{-1}x)\\
p_{\rho (g^*\circ-\id)}^2(y) &= -\rho^{2}p_{\rho^{-1}g}^2(-\rho^{-1}y)
\end{align*}
where $R_{\rho g}$ is defined in \eqref{eq:refl_prox} and $p_{\rho g}^2$ is defined in \eqref{eq:fcn_grad_rel_prox}.
\label{lem:refl_prox_conj_relations}
\end{lemma}

First, we show that the $z^k$ sequence in (primal) Douglas-Rachford \eqref{eq:DR} and the $v^k$ sequence in ADMM (i.e., dual Douglas-Rachford) in \eqref{eq:ADMM} differ by a factor only. This is well known \cite{EcksteinPhD}, but the relation is stated next with a simple proof.
\begin{proposition}
Assume that $\rho>0$ and $\gamma>0$ satisfy $\rho^{-1}=\gamma$. Further assume that $z^0 = \rho^{-1}v^0$. Then $z^k=\rho^{-1}v^{k}$ for all $k\geq 1$, where $\{z^k\}$ is the primal Douglas-Rachford sequence defined in \eqref{eq:DR} and the $\{v^k\}$ is the ADMM sequence is defined in \eqref{eq:ADMM}.
\end{proposition}
\begin{pf}
Lemma~\ref{lem:refl_prox_conj_relations} implies that
\begin{align*}
v^{k+1} &= (1-\alpha)v^k+\alpha R_{\rho (g^*\circ -\id)}R_{\rho f^*}v^k\\
&= (1-\alpha)v^k+\alpha \rho R_{\rho^{-1} g}(-\rho^{-1} (-\rho R_{\rho^{-1} f}(\rho^{-1}v^k)))\\
&= (1-\alpha)v^k+\alpha \rho R_{\rho^{-1} g}(R_{\rho^{-1} f}(\rho^{-1}v^k)))
\end{align*}
Multiply by $\rho^{-1}$, let $z^{k} = \rho^{-1}v^k$, and identify $\gamma=\rho^{-1}$ to get
\begin{align*}
z^{k+1} 
&= (1-\alpha)z^{k}+\alpha R_{\gamma g}(R_{\gamma f}(z^k))).
\end{align*}
This concludes the proof.
\end{pf}

There is also a tight relationship between the ADMM and Douglas-Rachford envelopes. Essentially, they have opposite signs.

\begin{proposition}
Assume that $\rho>0$ and $\gamma>0$ satisfy $\rho=\gamma^{-1}$ and that $z=\rho^{-1}v=\gamma v$. Then
\begin{align*}
F_{\rho}^{\rm{ADMM}}(v)&=-F_{\gamma}^{\rm{DR}}(z).
\end{align*}
\end{proposition}
\begin{pf}
Using Lemma~\ref{lem:refl_prox_conj_relations} several times, $\gamma=\rho^{-1}$, and $z=\rho^{-1}v$, we conclude that
\begin{align*}
F_{\rho}^{\rm{ADMM}}(v) &= (2\rho)^{-1}\left(\langle R_{\rho f^*}(v),v\rangle-p_{\rho f^*}^2(v)-p_{\rho(g^*\circ-\id)}^2(R_{\rho f^*}(v))\right)\\
&=(2\rho)^{-1}\Big(-\rho \langle R_{\rho^{-1}f}(\rho^{-1}v),v\rangle+\rho^2p_{\rho^{-1} (f\circ -\id)}(-\rho^{-1} v)\\
&\qquad\qquad\quad+\rho^2p_{\rho^{-1}g}(-\rho^{-1}(-\rho R_{\rho^{-1}f}(\rho^{-1}v)))\Big)\\
&=-\tfrac{\rho}{2}\left( \langle R_{\rho^{-1}f}(\rho^{-1}v),\rho^{-1} v\rangle-p_{\rho^{-1} f}^2(\rho^{-1} v)+p_{\rho^{-1}g}^2( R_{\rho^{-1}f}(\rho^{-1}v))\right)\\
&=-(2\gamma)^{-1}\left( \langle R_{\gamma f}(z),z\rangle-p_{\gamma f}^2(z)+p_{\gamma g}^2( R_{\gamma f}(z))\right)\\
&=-F_{\gamma}^{\rm{DR}}(z).
\end{align*}
This concludes the proof.
\end{pf}

This result implies that the ADMM envelope is concave when the DR envelope is convex, and vice versa. We know from Section~\ref{sec:DR} that the operator $S_1=R_{\rho f^*}$ is affine when $f^*$ is quadratic. This happens when
\begin{align*}
f(x)=\begin{cases} \tfrac{1}{2}\langle Hx,x\rangle+\langle h,x\rangle & {\hbox{if }} Ax=b\\
\infty & {\hbox{else}}
\end{cases}
\end{align*}
and $H$ is positive definite on the nullspace of $A$. From Proposition~\ref{prp:DR_bounds_general} and Proposition~\ref{prp:DR_bounds_id}, we conclude that, for an appropriate choice of $\rho$, the ADMM envelope is convex, which implies that the Douglas-Rachford envelope is concave.

\begin{remark}
The standard ADMM formulation is applied to solve problems of the form
\begin{align*}
\begin{tabular}{ll}
minimize & $\hat{f}(x)+\hat{g}(z)$\\
subject to & $Ax+Bz=c$
\end{tabular}
\end{align*}
Using infimal post-compositions, also called image functions, the dual of this is on the form \eqref{eq:ADMM_prob}, see e.g., \cite[Appendix~B]{gis_line_search} for details. So also this setting is implicitly considered.
\end{remark}

\section{The GAP Envelope}

\label{sec:GAP}

In this section, we provide an envelope function to a generalization of the classic alternating projections method in \cite{vonNeumann}. The generalization uses relaxed projections and is sometimes referred to as the method of alternating relaxed projections (MARP) \cite{MARP_Bauschke}, but we will refer to it as generalized alternating projections (GAP). The algorithm is analyzed in \cite{GAP_Gubin,GAP_Agmon,GAP_Motzkin,GAP_Eremin,GAP_Bregman} and a more general formulation is treated in \cite{Combettes2004}.

GAP solves feasibility problems with a finite number of nonempty closed and convex sets that have a nonempty intersection. Here, we consider feasibility problems with two sets:
\begin{align*}
{\hbox{find }} x\in C\cap D
\end{align*}
where $C\subset\reals^n$ and $D\subset\reals^n$ are nonempty closed and convex.

The generalized alternating projections method is given by
\begin{align}
x^{k+1}=(1-\alpha)x^k+\alpha P_{C}^{\alpha_2}P_{D}^{\alpha_1}x^k.
\label{eq:GAP}
\end{align}
where $P_{C}^{\alpha}$ is the relaxed projection in \eqref{eq:rel_proj}, and $\alpha\in(0,1]$ and $\alpha_1,\alpha_2\in(0,2]$. These assumptions imply that $P_{C}^{\alpha_2}$ is $\tfrac{\alpha_2}{2}$-averaged if $\alpha_2\in(0,2)$ and nonexpansive if $\alpha_2\in(0,2]$ (and similarly for $P_{D}^{\alpha_1}$). If $\alpha_1=2$ or $\alpha_2=2$, the composition $P_{C}^{\alpha_2}P_{D}^{\alpha_1}$ is nonexpansive and we need $\alpha\in(0,1)$ to arrive at an averaged iteration that guarantees convergence to a fixed-point. If $\alpha_1=\alpha_2=2$, the algorithm is Douglas-Rachford splitting (see Section~\ref{sec:DR}) applied to a feasibility problem. In this case, we have $\Pi_D({\rm{fix}}(P_{C}^{\alpha_2}P_{D}^{\alpha_1}))=C\cap D$. For all other feasible choices of $\alpha_1$ and $\alpha_2$, the fixed-point set satisfies ${\rm{fix}}(P_{C}^{\alpha_2}P_{D}^{\alpha_1})=C\cap D$. In either case, the algorithm performs an averaged iteration to find a fixed-point to the nonexpansive operator $P_{C}^{\alpha_2}P_{D}^{\alpha_1}$.

The algorithm is on the general form we consider and we identify $S_2$ in Assumption~\ref{ass:Si_prop} with $P_{C}^{\alpha_2}$ and $S_1$ with $P_D^{\alpha_1}$. We consider in particular the case when $S_1=P_{D}^{\alpha_1}$ is affine, i.e., $S_1=P(\cdot)+q$. This holds if $D$ is an affine set, i.e., if $D=\{x\in\reals^n~|~Ax=b\}$ for some linear operator $A$. Let $N$ denote the linear part of the projection onto the affine set $\Pi_D$, i.e., 
\begin{align}
N = \Pi_{D_0}
\label{eq:Ndef}
\end{align}
where $D_0=\{x\in\reals^n~|~Ax=0\}$, and let $d$ denote the constant part, to get $\Pi_Dx=Nx+d$. The operator $S_1$ then satisfies
\begin{align*}
S_1x&=P_{D}^{\alpha_1}x=(1-\alpha_1)x+\alpha_1\Pi_D=(1-\alpha_1)x+\alpha_1(Nx+d).
\end{align*}
This implies that $P$ and $q$ that define the affine operator $S_1=P(\cdot)+q$ satisfy
\begin{align}
\label{eq:P_GAP}P&=(1-\alpha_1)\id+\alpha_1N, & q&=\alpha_1d.
\end{align}

The GAP envelope function follows from the general envelope in \eqref{eq:env} and is given by
\begin{align*}
F_{\alpha_1,\alpha_2}^{\rm{GAP}}(x)=\tfrac{1}{2}\langle Px,x\rangle-p_{C}^{\alpha_2}(P_D^{\alpha_1}x)
\end{align*}
where $p_C^{\alpha_2}$ is defined in \eqref{eq:fcn_grad_rel_proj} and $P$ is from \eqref{eq:P_GAP}. Since $P_D^{\alpha_1}=Px+q$ and $\nabla p_{C}^{\alpha_2}=P_{C}^{\alpha_2}$, its gradient satisfies
\begin{align*}
\nabla F_{\alpha_1,\alpha_2}^{\rm{GAP}}(x)&= Px-P\nabla p_C^{\alpha_2}(Px+q)\\
&=P(x-P_C^{\alpha_2}P_D^{\alpha_1}x).
\end{align*}
So if $P$ is nonsingular, the stationary points of the GAP envelope coincides with the fixed-points of $P_C^{\alpha_2}P_D^{\alpha_1}$. The following proposition follows immediately from Proposition~\ref{prp:fp_sp_agree}.
\begin{proposition}
Suppose that $\alpha_1,\alpha_2\in(0,2]$ and that $\alpha_1\neq 1$. Then the set of stationary points to the gap envelope $F_{\alpha_1,\alpha_2}^{\rm{GAP}}$ is the fixed-point set of $P_{C}^{\alpha_2}P_{D}^{\alpha_1}$.
\label{prp:GAP_stat_points}
\end{proposition}

Next, we state some properties of the GAP envelope.
\begin{proposition}
Suppose that $\alpha_1\in(0,2]$ and $\alpha_2\in(0,2]$. Then the GAP envelope $F_{\alpha_1,\alpha_2}^{\rm{GAP}}$ satisfies
\begin{align*}
\tfrac{1}{2}\langle M(x-y),x-y\rangle&\leq F_{\alpha_1,\alpha_2}^{\rm{GAP}}(x)-F_{\alpha_1,\alpha_2}^{\rm{GAP}}(y)-\langle\nabla F_{\alpha_1,\alpha_2}^{\rm{GAP}}(y),x-y\rangle\\
&\leq \tfrac{1}{2}\langle L(x-y),x-y\rangle
\end{align*}
where
\begin{align}
M = \alpha_1(1-\alpha_1)(\id-N)
\label{eq:Mmat}
\end{align}
and
\begin{align}
L = (1-\alpha_1)(1+(\alpha_2-1)(1-\alpha_1))\id+\alpha_1(1+(\alpha_2-1)(2-\alpha_1))N
\label{eq:Lmat}
\end{align}
where $N$ is defined in \eqref{eq:Ndef}.
\label{prp:GAP_env_quad_bounds}
\end{proposition}
\begin{pf}
The operator $P_{C}^{\alpha_2}$ is $\tfrac{\alpha_2}{2}$-averaged and 1-negatively averaged (nonexpansive). So we can apply Theorem~\ref{thm:env_finer_prop} with $\delta_{\beta}=1$, $\delta_{\alpha}=\alpha_2-1$, and $P$ in \eqref{eq:P_GAP}. Using $N=N^2$ (which holds since $N$ is a projection onto a linear subspace), we conclude that
\begin{align*}
M&=P-P^2=(1-\alpha_1)\id+\alpha_1N-((1-\alpha_1)\id+\alpha_1N)^2\\
&=(1-\alpha_1)\id+\alpha_1N-((1-\alpha_1)^2\id+2\alpha_1(1-\alpha_1)N+\alpha_1^2N)\\
&=((1-\alpha_1)-(1-\alpha_1)^2)\id+(\alpha_1-(2\alpha_1-\alpha^2))N\\
&=((1-\alpha_1)-(1-2\alpha_1+\alpha_1^2))\id+(\alpha_1^2-\alpha_1))N\\
&=\alpha_1(1-\alpha_1)\id+\alpha_1(\alpha_1-1))N\\
&=\alpha_1(1-\alpha_1)(\id-N)
\end{align*}
and that
\begin{align*}
L&=P+(\alpha_2-1)P^2=(1-\alpha_1)\id+\alpha_1N+(\alpha_2-1)((1-\alpha_1)\id+\alpha_1N)^2\\
&=((1-\alpha_1)+(\alpha_2-1)(1-\alpha_1)^2)\id+(\alpha_1+(\alpha_2-1)(2\alpha_1(1-\alpha_1)+\alpha_1^2))N\\
&=(1-\alpha_1)(1+(\alpha_2-1)(1-\alpha_1))\id+\alpha_1(1+(\alpha_2-1)(2-\alpha_1))N.
\end{align*}
This concludes the proof.
\end{pf}

Since $N$ is a projection operator onto a linear subspace, it has only two distinct eigenvalues, namely zero and one. Therefore, there are only two distinct eigenvalues of $M$ and $L$ in \eqref{eq:Mmat} and \eqref{eq:Lmat}. Expressions for these eigenvalues are given in the following proposition.
\begin{proposition}
The eigenvalues of $M$ in \eqref{eq:Mmat} are 
\begin{align}
\lambda_{i}(M)=\begin{cases}
0 & {\hbox{for $i$ such that }} \lambda_{i}(N)=1\\
\alpha_1(1-\alpha_1)& {\hbox{for $i$ such that }} \lambda_{i}(N)=0\\
\end{cases}
\label{eq:Meig}
\end{align}
and the eigenvalues of $L$ in \eqref{eq:Lmat} are
\begin{align}
\lambda_{i}(L)=\begin{cases}
\alpha_2 & {\hbox{for $i$ such that }} \lambda_{i}(N)=1\\
(1-\alpha_1)(1+(\alpha_2-1)(1-\alpha_1))& {\hbox{for $i$ such that }} \lambda_{i}(N)=0\\
\end{cases}
\label{eq:Leig}
\end{align}
with $N$ defined in \eqref{eq:Ndef}.
\label{prp:LM_eig}
\end{proposition}
\begin{pf}
First note that $\lambda_{i}(a_1\id+a_2N)=a_1+a_2\lambda_i(N)$. This implies that $\lambda_{i}(M)=\alpha_1(1-\alpha_1)(1-\lambda_{i}(N))$, and \eqref{eq:Meig} is proven. It also implies that
\begin{align*}
\lambda_{i}(L) = (1-\alpha_1)(1+(\alpha_2-1)(1-\alpha_1))+\alpha_1(1+(\alpha_2-1)(2-\alpha_1))\lambda_{i}(N).
\end{align*}
For $\lambda_{i}(N)=0$, we see that \eqref{eq:Leig} holds. In the case of $\lambda_i(N)=1$, we conclude that
\begin{align*}
\lambda_{i}(L)&=(1-\alpha_1)(1+(\alpha_2-1)(1-\alpha_1))+\alpha_1(1+(\alpha_2-1)(2-\alpha_1))\\
&=1-\alpha_1+\alpha_2(1-\alpha_1)^2-(1-\alpha_1)^2+\alpha_1+\alpha_1\alpha_2(2-\alpha_1)-\alpha_1(2-\alpha_1)\\
&=1+\alpha_2(1-2\alpha_1+\alpha_1^2)-1+2\alpha_1-\alpha_1^2+\alpha_1\alpha_2(2-\alpha_1)-2\alpha_1-\alpha_1^2\\
&=\alpha_2(1-2\alpha_1+\alpha_1^2)+\alpha_2(2\alpha_1-\alpha_1^2)\\
&=\alpha_2.
\end{align*}
This concludes the proof.
\end{pf}

Using this, we can show that for $\alpha_1\in[1,2]$, the GAP envelope is convex on the nullspace of $A$ and concave on its orthogonal complement, the rangespace of $A^*$.
\begin{proposition}
  Let $\mathcal{N}(A)$ denote the nullspace of $A$ and let $\mathcal{R}(A^*)$ denote its orthogonal complement, the rangespace of $A^*$. Then the GAP envelope is convex and $\alpha_2$-smooth when restricted to $\mathcal{R}(A^*)$. If $\alpha_1\in[1,2]$, the GAP envelope is concave and $\alpha_1(\alpha_1-1)$-smooth when restricted to $\mathcal{N}(A)$.
\end{proposition}
\begin{pf}
The subspace $\mathcal{R}(A^*)$ is spanned by the eigenvectors corresponding to $\lambda_i(N)=1$. Therefore, Proposition~\ref{prp:LM_eig} implies that for all $x,y\in\mathcal{R}(A^*)$, the lower bound in Proposition~\ref{prp:GAP_env_quad_bounds} becomes $\langle M(x-y),x-y\rangle=0$ and the upper bound in Proposition~\ref{prp:GAP_env_quad_bounds} satisfies $\langle L(x-y),x-y\rangle=\alpha_2\|x-y\|^2$. This proves the first claim.

The second claim is proven similarly. The subspace $\mathcal{N}(A)$ is spanned by the eigenvectors corresponding to $\lambda_{i}(N)=0$. Therefore, Proposition~\ref{prp:LM_eig} implies that for all $x,y\in\mathcal{N}(A)$, the lower bound in Proposition~\ref{prp:GAP_env_quad_bounds} becomes $\langle M(x-y),x-y\rangle=\alpha_1(1-\alpha_1)\|x-y\|^2$ and the upper bound in Proposition~\ref{prp:GAP_env_quad_bounds} satisfies $\langle L(x-y),x-y\rangle=(1-\alpha_1)(1+(\alpha_2-1)(1-\alpha_1))\|x-y\|^2$. Noting that $(1-\alpha_1)(1+(\alpha_2-1)(1-\alpha_1))\leq 0$ when $\alpha_1\in[1,2]$ and $\alpha_2\in(0,2]$ proves the second claim. 
\end{pf}

The following proposition is a straightforward consequence of Proposition~\ref{prp:GAP_env_quad_bounds} and Proposition~\ref{prp:LM_eig} and is stated without a proof.
\begin{proposition}
Suppose that $\alpha_1\in(0,2]$ and $\alpha_2\in(0,2]$. Then the GAP envelope $F_{\alpha_1,\alpha_2}^{\rm{GAP}}$ satisfies
\begin{align*}
\tfrac{\beta_l}{2}\|x-y\|^2&\leq F_{\alpha_1,\alpha_2}^{\rm{GAP}}(x)-F_{\alpha_1,\alpha_2}^{\rm{GAP}}(y)-\langle\nabla F_{\alpha_1,\alpha_2}^{\rm{GAP}}(y),x-y\rangle\leq \tfrac{\beta_u}{2}\|x-y\|^2
\end{align*}
where $\beta_l=\min((1-\alpha_1)\alpha_1,0)$ and $\beta_u=\max((1-\alpha_1)(1+(\alpha_2-1)(1-\alpha_1)),\alpha_2)$. If in addition $\alpha_1\in(0,1]$, then it is convex.
\label{prp:GAPenv_bounds_id}
\end{proposition}

If the first relaxed projection is under-relaxed, i.e., if $\alpha_1\in(0,1]$, then the GAP envelope is convex. From Proposition~\ref{prp:GAP_stat_points}, we also know that if $\alpha_1\neq 1$ its set of stationary points is the fixed-point set of $P_{C}^{\alpha_2}P_{D}^{\alpha_1}$. For convex functions, all stationary points are minimizers. This therefore implies that all convex feasibility problems where one set is affine, can be solved by minimizing the smooth convex GAP envelope function by setting $\alpha_1\in(0,1)$. In Section~\ref{sec:cone}, we will see that most convex optimization problems can actually be cast on this feasibility form.

\section{Conclusions}

We have presented a unified framework for envelope functions. Special cases include the Moreau envelope, the forward-backward envelope, the Douglas-Rachford envelope, and the ADMM envelope. We also presented a new envelope function, namely the generalized alternating projections (GAP) envelope. Under additional assumptions, we have provided quadratic upper and lower bounds to the general envelope function. These coincide with or sharpen corresponding results for the known special cases in the literature. 

\section{Acknowledgments}

Both authors are financially supported by the Swedish Foundation for Strategic Research and members of the LCCC Linneaus Center at Lund University.

\bibliographystyle{plain}
\bibliography{references}

\appendix

\section{Proof to Theorem~\ref{thm:env_finer_prop}}
\label{app:env_finer_prop_pf}
First, we establish that
\begin{align}
-\delta_{\alpha}\|x-y\|_{P^2}^2\leq \langle P\nabla f_2(Px+q)-P\nabla f_2(Py+q),x-y\rangle\leq \delta_{\beta} \|x-y\|_{P^2}^2.
\label{eq:comp_neg_mono}
\end{align}
We have
\begin{align*}
\langle P\nabla f_2(Px+q)-P&\nabla f_2(Py+q),x-y\rangle\\
&=\langle \nabla f_2(Px+q)-\nabla f_2(Py+q),P(x-y)\rangle\\
&=\langle \nabla f_2(Px+q)-\nabla f_2(Py+q),(Px+q)-(Py+q))\rangle
\end{align*}
This implies that
\begin{align*}
-(2\alpha-1)\|x-y\|_{P^2}^2&=-(2\alpha-1)\|(Px+q)-(Py-q)\|^2\\
&\leq\langle P\nabla f_2(Px+q)-P\nabla f_2(Py+q),x-y\rangle\\
&\leq(2\beta-1)\|(Px+q)-(Py-q)\|^2\\
&=(2\beta-1)\|x-y\|_{P^2}^2
\end{align*}
where Lemma~\ref{lem:id_sub_grad_mono_avg} and Lemma~\ref{lem:id_sub_grad_mono_negavg} are used in the inequalities. Recalling that $\delta{\alpha}=2\alpha-1$ and $\delta_{\beta}=2\beta-1$, this shows that \eqref{eq:comp_neg_mono} holds. Further, for any $\delta\in\reals$ we have 
\begin{align}
\nonumber \langle \nabla F(x)-\nabla F(y),x-y\rangle &=\langle P(x-\nabla f_2\nabla f_1(x))-P(x-\nabla f_2\nabla f_1(y)),x-y\rangle \\
\nonumber &=\langle P(x-y),x-y\rangle\\
\nonumber &\quad-\langle P\nabla f_2(Px+q)-P\nabla f_2(Py+q),x-y\rangle\\
\nonumber &=\langle(P-\delta P^2)(x-y),x-y\rangle+\delta\|x-y\|_{P^2}^2\\
\label{eq:grad_sp_eq}&\quad-\langle P\nabla f_2(Px+q)-P\nabla f_2(Py+q),x-y\rangle.
\end{align}
Let $\delta=-\delta_{\alpha}$, then \eqref{eq:grad_sp_eq} and \eqref{eq:comp_neg_mono} imply
\begin{align*}
\langle \nabla F(x)-\nabla F(y),x-y\rangle &\leq \langle(P+\delta_{\alpha} P^2)(x-y),x-y\rangle.
\end{align*}
Let $\delta=\delta_{\beta}$, then \eqref{eq:grad_sp_eq} and \eqref{eq:comp_neg_mono} imply
\begin{align*}
\langle \nabla F(x)-\nabla F(y),x-y\rangle &\geq \langle(P-\delta_{\beta} P^2)(x-y),x-y\rangle.
\end{align*}
Applying Lemma~\ref{lem:f_bounds} in Appendix~\ref{app:lemmas} gives the result.

\section{Proof to Lemma~\ref{lem:refl_prox_conj_relations}}
\label{app:refl_prox_conj_relations}

Using the Moreau decomposition \cite[Theorem~14.3]{bauschkeCVXanal}
\begin{align*}
\prox_{\rho g^*}(x) = x-\rho\prox_{\rho^{-1}g}(\rho^{-1}x),
\end{align*}
we conclude that 
\begin{align*}
R_{\rho g^*}(x) &= 2\prox_{\rho g^*}(x)-x\\
&=2(x-\rho\prox_{\rho^{-1}g}(\rho^{-1}x))-x\\
&=-\rho\left(2(\prox_{\rho^{-1}g}(\rho^{-1}x))-(\rho^{-1}x)\right)\\
&=-\rho R_{\rho^{-1}g}(\rho^{-1}x)
\end{align*}
and
\begin{align*}
R_{\rho (g^*\circ -\id)}(x) &= 2\prox_{\rho(g^*\circ-\id)}(x)-x\\
&=-2\prox_{\rho g^*}(-x)-x\\
&=-2(-x-\rho \prox_{\rho^{-1}g}(-\rho^{-1} x))-x\\
&=2\rho \prox_{\rho^{-1}g}(-\rho^{-1} x))+x\\
&=\rho(2 \prox_{\rho^{-1}g}(-\rho^{-1} x)-(-\rho^{-1}x))\\
&=\rho R_{\rho^{-1} g}(-\rho^{-1}x).
\end{align*}
To show the third claim, we first derive an expression for $r_{\rho(g^*\circ -\id)}^*$. We have
\begin{align*}
r_{\rho (g^*\circ-\id)}^*(y) &= (\rho (g^*\circ-\id)+\tfrac{1}{2}\|\cdot\|^2)^*(y)\\
&=\sup_{z}\{\langle y,z\rangle-\rho \sup_{x}\{\langle z,x\rangle-g(-x)\}-\tfrac{1}{2}\|z\|^2\}\\
&=\sup_{z}\{\langle y,z\rangle+\rho \inf_{x}\{\langle z,-x\rangle+g(-x)\}-\tfrac{1}{2}\|z\|^2\}\\
&=\sup_{z}\{\langle y,z\rangle+\rho \inf_{v}\{\langle z,v\rangle+g(v)\}-\tfrac{1}{2}\|z\|^2\}\\
&=\sup_{z}\inf_v\{\langle y,z\rangle+\rho \langle z,v\rangle+\rho g(v)-\tfrac{1}{2}\|z\|^2\}\\
&=\inf_v\sup_{z}\{\langle y+\rho v,z\rangle+\rho g(v)-\tfrac{1}{2}\|z\|^2\}\\
&=\inf_v\{\tfrac{1}{2}\| y+\rho v\|^2+\rho g(v)\}\\
&=\inf_v\{\langle y,\rho v\rangle+\tfrac{1}{2}\| \rho v\|^2+\rho g(v)\}+\tfrac{1}{2}\|y\|^2\\
&=-\sup_v\{\langle -y,\rho v\rangle-\tfrac{1}{2}\| \rho v\|^2-\rho g(v)\}+\tfrac{1}{2}\|y\|^2\\
&=-\rho^2\sup_v\{\langle -\rho^{-1}y,v\rangle-\tfrac{1}{2}\| v\|^2-\rho^{-1} g(v)\}+\tfrac{1}{2}\|y\|^2\\
&=-\rho^2r_{\rho^{-1}g}^*(-\rho^{-1}y)+\tfrac{1}{2}\|y\|^2,
\end{align*}
where the sup-inf swap is valid by the minimax theorem in \cite{minimax_Sion} since we can construct a compact set for the $z$ variable due to strong convexity of $\|\cdot\|^2$. This implies that
\begin{align*}
p_{\rho (g^*\circ-\id)}^2(y) &= 2r_{\rho (g^*\circ-\id)}^*(y)-\tfrac{1}{2}\|y\|^2\\
&=-2\rho^{2}r_{\rho^{-1}g}^*(-\rho^{-1}y)+\tfrac{1}{2}\|y\|^2\\
&=-\rho^{2}(2r_{\rho^{-1}g}^*(-\rho^{-1}y)-\tfrac{1}{2}\|-\rho^{-1}y\|^2)\\
&=-\rho^{2}p_{\rho^{-1}g}^2(-\rho^{-1}y).
\end{align*}
This concludes the proof.

\section{Technical Lemmas}
\label{app:lemmas}

\begin{lemma}
Assume that $f~:~\reals^n\to\reals$ is differentiable and that $M~:~\reals^n\to\reals^n$ and $L~:~\reals^n\to\reals^n$ are linear operators. Then
\begin{align}
-\tfrac{1}{2}\langle M(x-y),x-y\rangle\leq f(x)-f(y)-\langle\nabla f(y),x-y\rangle\leq\tfrac{1}{2}\langle L(x-y),x-y\rangle
\label{eq:quad_bound}
\end{align}
if and only if
\begin{align}
-\langle M(x-y),x-y\rangle\leq\langle\nabla f(x)-\nabla f(y),x-y\rangle\leq\langle L(x-y),x-y\rangle
\label{eq:grad_bound}
\end{align}
\label{lem:f_bounds}
\end{lemma}
\begin{pf}
Adding two copies of \eqref{eq:quad_bound} with $x$ and $y$ interchanged gives
\begin{align}
-\langle M(x-y),x-y\rangle\leq \langle\nabla f(x)-f(y),x-y\rangle\leq\langle L(x-y),x-y\rangle.
\end{align}
This shows that \eqref{eq:quad_bound} implies \eqref{eq:grad_bound}. To show the other direction, we use integration. Let $h(\tau)=f(x+\tau(y-x))$, then
\begin{align*}
\nabla h(\tau) = \langle y-x,\nabla f(x+\tau(y-x))\rangle
\end{align*}
 since $f(y)=h(1)$ and $f(x)=h(0)$, we get
\begin{align*}
f(y)-f(x)&=h(1)-h(0)=\int_{0}^{1}\nabla h(\tau)d\tau
=\int_{0}^1\langle y-x,\nabla f(x+\tau(y-x))\rangle d\tau
\end{align*}
Therefore
\begin{align*}
f(y)-f(x)-\langle \nabla f(x),y-x\rangle
&=\int_0^1\langle \nabla f(x+\tau(y-x)),y-x\rangle d\tau-\langle \nabla f(x),y-x\rangle\\
&=\int_0^1\langle \nabla f(x+\tau(y-x))-\nabla f(x),y-x\rangle d\tau\\
&=\int_0^1\tau^{-1}\langle \nabla f(x+\tau(y-x))-\nabla f(x),\tau(y-x)\rangle d\tau\\
&=\int_0^1\tau^{-1}\langle \nabla f(x+\tau(y-x))-\nabla f(x),(x+\tau(y-x))-x\rangle d\tau.
\end{align*}
Using the upper bound in \eqref{eq:grad_bound}, we get
\begin{align*}
\int_0^1\tau^{-1}\langle \nabla f(x+\tau(y-x))&-\nabla f(x),(x+\tau(y-x))-x\rangle d\tau\\
&\leq \int_0^1\tau^{-1}\langle  L\tau (x-y),\tau(x-y)\rangle d\tau\\
&=\langle L (x-y),x-y\rangle\int_0^1\tau d\tau\\
&=\tfrac{1}{2}\langle L (x-y),x-y\rangle.
\end{align*}
Similarly, using the lower bound in \eqref{eq:grad_bound}, we get
\begin{align*}
\int_0^1\tau^{-1}\langle \nabla f(x+\tau(y-x))&-\nabla f(x),(x+\tau(y-x))-x\rangle d\tau\\
&\geq -\int_0^1\tau^{-1}\langle  M\tau (x-y),\tau(x-y)\rangle d\tau\\
&=-\langle M(x-y),x-y\rangle\int_0^1\tau d\tau\\
&=-\tfrac{1}{2}\langle M (x-y),x-y\rangle.
\end{align*}
This concludes the proof.
\end{pf}

\begin{lemma}
Assume that $f~:~\reals^n\to\reals$ is differentiable and that $L$ is positive definite. Then that $f$ is $L$-smooth, i.e., that $f$ satisfies
\begin{align}
|f(x)-f(y)-\langle\nabla f(y),x-y\rangle|\leq \tfrac{\beta}{2}\|x-y\|_L^2
\label{eq:quad_bound_scaled_equal}
\end{align}
holds for all $x,y\in\reals^n$ is equivalent to that $\nabla f$ is $\beta$-Lipschitz continuous w.r.t. $\|\cdot\|_L$, i.e., that
\begin{align}
\|\nabla f(x)-\nabla f(y)\|_{L^{-1}}\leq \beta\|x-y\|_{L}
\label{eq:scaled_Lipschitz}
\end{align}
holds for all $x,y\in\reals^n$.
\label{lem:quad_bound_Lipschitz}
\end{lemma}
\begin{pf}
We start by proving the result using the induced norm $\|\cdot\|$ only, i.e., in the Hilbert space setting. (This covers, e.g., the setting with inner-product $\langle x,y\rangle_H=\langle Hx,y\rangle$ and scaled norm $\|\cdot\|_H=\sqrt{\langle x,y\rangle_H}$ that will be used later.) To do this, we introduce the functions $h:=\tfrac{1}{\beta} f$ and $r:=\tfrac{1}{2}(h+\tfrac{1}{2}\|\cdot\|^2)$.

Since $L=\id$ in the norm, the condition \eqref{eq:scaled_Lipschitz} is $\beta$-Lipschitz continuity of $\nabla f$ (w.r.t. $\|\cdot\|$). This is equivalent to that $\nabla h=\tfrac{1}{\beta}\nabla f$ is nonexpansive, which by \cite[Proposition~4.2]{bauschkeCVXanal} is equivalent to that $\tfrac{1}{2}(\nabla h+\id)=\nabla \left(\tfrac{1}{2}(h+\tfrac{1}{2}\|\cdot\|^2)\right)=\nabla r$ is firmly nonexpansive (or equivalently 1-cocoercive). This, is equivalent to (see \cite[Theorem~2.1.5]{NesterovLectures} and \cite[Definition~4.4]{bauschkeCVXanal}) that:
\begin{align*}
0\leq r(x)-r(y)-\langle \nabla r(y),x-y\rangle\leq \tfrac{1}{2}\|x-y\|^2.
\end{align*}
holds for all $x,y\in\reals^n$. Multiplying by 2 and using $2r=h+\tfrac{1}{2}\|\cdot\|^2$, this is equivalent to that
\begin{align*}
  0&\leq h(x)-h(y)-\langle\nabla h(y),x-y\rangle+\tfrac{1}{2}(\|x\|^2-\|y\|^2-2\langle y,x-y\rangle)\\
&=h(x)-h(y)-\langle\nabla h(y),x-y\rangle+\tfrac{1}{2}\|x-y\|^2\leq \|x-y\|^2.
\end{align*}
Multiplying by $\beta$ and using $f=\beta h$, this is equivalent to
\begin{align*}
  -\tfrac{\beta}{2}\|x-y\|&\leq f(x)-f(y)-\langle\nabla f(y),x-y\rangle\leq \tfrac{\beta}{2}\|x-y\|^2.
\end{align*}
This chain of equivalences show that the conditions are equivalent when $L=\id$.

Next, we show that the scaled version holds. To do this, introduce the space $\mathbb{H}_H$ with inner-product $\langle x,y\rangle_H=\langle Hx,y\rangle$ and induced norm $\|\cdot\|_H=\sqrt{\langle Hx,x\rangle}$ and the space  $\mathbb{E}_L$ inner-product $\langle x,y\rangle$ and induced norm $\|\cdot\|_L=\sqrt{\langle Lx,x\rangle}$. Further let $H=L$ and define $f_h~:~\mathbb{H}_H\to\reals$ and $f_l~:~\mathbb{E}_L\to\reals$ that satisfy $f_h(x)=f_l(x)$ for all $x\in\reals^n$. We have already shown that \eqref{eq:quad_bound_scaled_equal} and \eqref{eq:scaled_Lipschitz} are equivalent for $f_h$ that is defined on the Hilbert space $\mathbb{H}_H$. To show that it also holds for $f_l$ defined on $\mathbb{E}_L$, we show that the conditions \eqref{eq:quad_bound_scaled_equal} and \eqref{eq:scaled_Lipschitz} are equivalent if defined for $f_h$ on $\mathbb{H}_H$ and if defined for $f_l$ on $\mathbb{E}_L$, when $L=H$.

By definition of the gradient, $\nabla f_l$ and $\nabla f_h$ must satisfy
\begin{align*}
\langle\nabla f_l(y),x-y\rangle = \langle\nabla f_h(y),x-y\rangle_H=\langle H\nabla f_h(y),x-y\rangle
\end{align*}
for all $x,y\in\reals^n$. This implies that $\nabla f_h=H^{-1}\nabla f_l=L^{-1}\nabla f_l$. Therefore that \eqref{eq:quad_bound_scaled_equal} holds for $f_l$ on $\mathbb{E}_L$ is equivalent to that it holds for $f_h$ on $\mathbb{H}_H$.

Further,
\begin{align*}
\|\nabla f_h(x)-\nabla f_h(y)\|_H^2&=\langle \nabla f_h(x)-\nabla f_h(y),\nabla f_h(x)-\nabla f_h(y)\rangle_H\\
&=\langle L^{-1}(\nabla f(x)-\nabla f(y)),L^{-1}(\nabla f(x)-\nabla f(y))\rangle_L\\
&=\langle \nabla f(x)-\nabla f(y),\nabla f(x)-\nabla f(y)\rangle_{L^{-1}}\\
&=\|\nabla f(x)-\nabla f(y)\|_{L^{-1}}^2.
\end{align*}
So that \eqref{eq:scaled_Lipschitz} holds for $f_l$ on $\mathbb{E}_L$ is equivalent to that it holds for $f_h$ on $\mathbb{H}_H$. This concludes the proof.
\end{pf}

\begin{lemma}
Assume that $f$ is differentiable. Then $\nabla f$ is $\alpha$-averaged with $\alpha\in(0,1]$ if and only if
\begin{align}
-(2\alpha-1)\|x-y\|^2\leq\langle \nabla f(x)-\nabla f(y),x-y\rangle\leq\|x-y\|^2.
\label{eq:grad_prop_avg}
\end{align}
\label{lem:id_sub_grad_mono_avg}
\end{lemma}
\begin{pf}
The operator $\nabla f$ is $\alpha$-averaged if and only if $\nabla f=(1-\alpha)\id+\alpha R$ for some nonexpansive operator $R$. Therefore, $\nabla f$ is $\alpha$-averaged if and only if $\nabla f-(1-\alpha)\id$ is $\alpha$-Lipschitz continuous, since $\nabla f-(1-\alpha)\id=\alpha R$. Letting $g := f-\tfrac{1-\alpha}{2}\|\cdot\|^2$, we get $\nabla g=\alpha R$. Therefore $\nabla g$ is $\alpha$-Lipschitz. According to Lemma~\ref{lem:quad_bound_Lipschitz} this is equivalent to that
\begin{align*}
|g(x)-g(y)-\langle\nabla g(y),x-y\rangle|\leq\tfrac{\alpha}{2}\|x-y\|^2
\end{align*}
or equivalently
\begin{align*}
|f(x)-f(y)-\langle\nabla f(y),x-y\rangle-\tfrac{1-\alpha}{2}\|x-y\|^2|\leq\tfrac{\alpha}{2}\|x-y\|^2
\end{align*}
which is equivalent to
\begin{align*}
-\tfrac{2\alpha-1}{2}\|x-y\|^2\leq f(x)-f(y)-\langle\nabla f(y),x-y\rangle\leq\tfrac{1}{2}\|x-y\|^2.
\end{align*}
Applying Lemma~\ref{lem:f_bounds} gives the result.
\end{pf}
\begin{lemma}
Assume that $f$ is differentiable. Then $\nabla f$ is $\beta$-negatively averaged with $\beta\in(0,1]$ if and only if
\begin{align}
-\|x-y\|^2\leq\langle \nabla f(x)-\nabla f(y),x-y\rangle\leq(2\beta-1)\|x-y\|^2.
\label{eq:grad_prop_negavg}
\end{align}
\label{lem:id_sub_grad_mono_negavg}
\end{lemma}
\begin{pf}
This follows immediately from \ref{lem:id_sub_grad_mono_avg} since $-\nabla f$ is $\beta$-averaged by definition.
\end{pf}

\begin{lemma}
Suppose that $P$ is a linear self-adjoint and nonexpansive operator with largest eigenvalue $\lambda_{\max}(P)=L$ and smallest eigenvalue $\lambda_{\min}(P)=m$, satisfying $-1\leq m\leq L\leq 1$. Further suppose that $\delta\in[-1,1]$ and let $j$ be the index that minimizes $|\tfrac{1}{2\delta}-\lambda_i(P)|$, i.e., $j=\argmin_i(|\tfrac{1}{2\delta}-\lambda_i(P)|)$. The smallest eigenvalue of $P-\delta P^2$ satisfies the following: 
\begin{enumerate}[(i)]
\item if $\delta\in[0,1]$, then $\lambda_{\min}(P-\delta P^2)=\min(m-\delta m^2,L-\delta L^2)$ 
\item if $\delta\in[-0.5,0]$, then $\lambda_{\min}(P-\delta P^2)=m-\delta m^2$ 
\item  if $\delta\in[-1,-0.5]$, then $\lambda_{\min}(P-\delta P^2)=\lambda_j(P)-\delta\lambda_j(P)^2$
\end{enumerate}
\label{lem:smallest_eigenvalue}
\end{lemma}
\begin{pf}
From the spectral theorem it follows that the eigenvalues to $\lambda_i(P-\delta P^2) = \lambda_i(P)-\delta\lambda_i(P)^2$. So we need to find the $\lambda_i(P)$ that minimizes the function $\psi(\lambda)=\lambda-\delta\lambda^2$, where $\lambda_i(P)\in[-1,1]$ for different $\delta$.

For $\delta\in[0,1]$, the function $\psi$ is concave, and the minimum is found in either of the end points, so $\lambda_{\min}(P-\delta P^2)=\min(m-\delta m^2,L-\delta L^2)$. This shows {\it{(i)}}. If instead $\delta\in[-1,0)$ the function $\psi$ is convex. The unconstrained minimum is at $\tfrac{1}{2\delta}$. Then, since the level sets of $\psi$ are symmetric around $\tfrac{1}{2\delta}$, the constrained minimum is the eigenvalue $\lambda_i(P)$ closest to $\tfrac{1}{2\delta}$. For $\delta\in[-0.5,0)$ this is $\lambda_{\min}(P)=m$, and for $\delta\in[-1,-0.5]$ this is $\lambda_j(P)$. This concludes the proof.
\end{pf}

\begin{lemma}
Suppose that $P$ is a linear self-adjoint and nonexpansive operator with largest eigenvalue $\lambda_{\max}(P)=L$ and smallest eigenvalue $\lambda_{\min}(P)=m$, satisfying $-1\leq m\leq L\leq 1$. Further suppose that $\delta\in[-1,1]$ and let $j$ be the index that minimizes $|\tfrac{1}{2\delta}+\lambda_i(P)|$, i.e., $j=\argmin_i(|\tfrac{1}{2\delta}+\lambda_i(P)|)$. The largest eigenvalue of $P+\delta P^2$ satisfies the following:
\begin{enumerate}[(li)]
\item if $\delta\in[-0.5,1]$, then $\lambda_{\max}(P+\delta P^2)=L+\delta L^2$ 
\item  if $\delta\in[-1,-0.5]$, then $\lambda_{\max}(P+\delta P^2)=\lambda_j(P)+\delta\lambda_j(P)^2$
\end{enumerate}
\label{lem:largest_eigenvalue}
\end{lemma}
\begin{pf}
From the spectral theorem it follows that the eigenvalues to $\lambda_i(P+\delta P^2) = \lambda_i(P)+\delta\lambda_i(P)^2$. So we need to find the $\lambda_i(P)$ that maximizes the function $\psi(\lambda)=\lambda+\delta\lambda^2$, where $\lambda_i(P)\in[-1,1]$ for different $\delta$.

For $\delta\in[0,1]$, the function $\psi$ is convex, and the maximum is found in either of the end points. The function $\psi$ is monotonically increasing on $[-1,1]$, so the maximum is found at $L+\delta L^2$. For $\delta\in[-1,0)$, the function $\psi$ is concave. Its unconstrained maximum is at $\tfrac{1}{-2\delta}$. Since the level sets of $\psi$ are symmetric around $\tfrac{1}{-2\delta}$, the constrained maximum is the eigenvalue closest to $\tfrac{1}{-2\delta}$. For $\delta\in[-0.5,0)$, this is $\lambda_{\max}(P)=L$, and for $\delta\in[-1,-0.5]$ this is $\lambda_j(P)$. This concludes the proof.
\end{pf}


\end{document}